\documentclass[final,twocolumn,twoside]{IEEEtran}

\usepackage{amsmath,amssymb,theorem}
\usepackage{graphicx,floatflt}
\usepackage{algorithm,algorithmic}
\usepackage{psfrag}
	
\usepackage{color,xspace}
\usepackage{subfigure}

\newtheorem{theorem}{Theorem}[section]
\newtheorem{lemma}[theorem]{Lemma}

\newtheorem{remark}[theorem]{Remark}

\newtheorem{definition}[theorem]{Definition}
\newtheorem{problem}[theorem]{Problem}

\newcommand{\real}{{\mathbb{R}}}

\newcommand{\realnonnegative}{\mathbb{R}_{\ge 0}}
\newcommand{\integernonnegative}{\mathbb{Z}_{\ge 0}}

\newcommand{\NN}{{\mathcal{N}}}

\renewcommand{\natural}{{\mathbb{N}}}

\renewcommand{\tilde}{\widetilde}

\renewcommand{\epsilon}{\varepsilon}

\newcommand{\until}[1]{\{1,\dots, #1\}}

\renewcommand{\hat}{\widehat}

\newcommand{\Nouti}{\NN_i^\text{out}}

\newcommand{\clovar}{\chi}

\newcommand{\sumvar}{\phi}
\newcommand{\solvar}{\varphi}

\newcommand{\oprocendsymbol}{\hbox{$\bullet$}}
\newcommand{\oprocend}{\relax\ifmmode\else\unskip\hfill\fi\oprocendsymbol}

\newcommand{\longthmtitle}[1]{\mbox{}\textup{\textbf{(#1)}}}

\begin{document}

\title{Robust Dynamic Event-Triggered Coordination \\ With a Designable Minimum Inter-Event Time}

\author{James Berneburg \qquad Cameron Nowzari \thanks{The authors are
    with the Department of Electrical and Computer Engineering,
    George Mason University, Fairfax, VA 22030, USA, {\tt\small
      \{jbernebu,cnowzari\}}@gmu.edu}}
\maketitle

\begin{abstract}
This paper revisits the classical multi-agent average consensus problem for which many different event-triggered control strategies have been proposed over the last decade. Many of the earliest versions of these works conclude asymptotic stability without proving that Zeno behavior, or deadlocks, do not occur along the trajectories of the system. More recent works
that resolve this issue either: (i) propose the use of a dwell-time that forces inter-event times to be lower-bounded away from zero but sacrifice asymptotic convergence in exchange for practical convergence (or convergence to a neighborhood); (ii) guarantee non-Zeno behaviors and asymptotic convergence but do not provide a positive minimum inter-event time guarantee; or (iii) are not fully distributed. 
Additionally, the overwhelming majority of these works provide no form of robustness analysis on the event-triggering strategy. More specifically, if arbitrarily small disturbances can remove the non-Zeno property then the theoretically correct algorithm may not actually be implementable. Instead, this work for the first time presents a fully distributed, robust, dynamic event-triggered algorithm, for general directed communication networks, for which a desired positive minimum inter-event time can be chosen by each agent in a distributed fashion. Simulations illustrate our results.
\end{abstract} 

\section{Introduction}


Systems composed of individually controlled agents are increasingly common and a very active area of research. Such systems are designed for many different applications including the coordination of unmanned air vehicles, distributed reconfigurable sensor networks, and attitude alignment for satellites, etc; see~\cite{RosRmm2004} and~\cite{WrRwbEma2007} and their references. These are often intended to fulfill some coordinated task, but require distributed control to be scaled with large systems. In this case, communication limitations, such as wireless bandwidth, mean that agents cannot be assumed to have continuous access to others' states.
Therefore, many works have recently considered communication to be a limited resource, where individual agents must autonomously schedule when to take various actions, rather than doing so periodically or continuously.

A common solution to these types of problems comes in the form of event-triggered coordination, where actions occur at specific instances of time when some event condition is satisfied, such as when an error state~\cite{XyKlDvdKhj17} or a clock state~\cite{CdpRp2017} hits some threshold. A similar strategy is self-triggered control, where the controller uses state information to schedule events ahead of time. An introduction to these ideas for single-plant systems is found in~\cite{WpmhhKhjPt2012}. 

One potential problem in event-triggered coordination is the Zeno phenomenon, where the number of events triggered goes to infinity in a finite time period. This is problematic as it asks for solutions that cannot be realized by actual devices. A way to prevent this problem is to design triggering conditions that guarantee a positive minimum inter-event time (MIET) exists. Note that this is different from first designing an event-triggering condition, and then afterwards \emph{forcing} a minimum inter-event time (dwell-time), which has drawbacks as we discuss later. For example, the self-triggered strategy in~\cite{CdpPf2013} enforces a MIET and so only guarantees convergence to a set.

As noted in~\cite{WpmhhKhjPt2012}, the existence of a positive MIET is important to ensure that the event-triggering mechanism does not become unimplementable because it requires actions to be taken arbitrarily fast. This issue has been addressed recently for single-plant systems, e.g., in~\cite{AvpMmj2018}, where a general method for achieving stabilization based on Lyapunov functions is developed, and~\cite{RpPtDnAa2015}, where a general framework for event-triggering mechanisms is provided using the hybrid systems formalism of~\cite{RgRgsArt2012}. 
{Hybrid systems formulations seem especially useful in networked control systems because they conveniently describe systems with continuous-time dynamics and discrete-time memory updates via communication.} 
However, this is still a major challenge for multi-agent systems with distributed information.

To address this, we turn to a simple but widely applicable canonical problem: multi-agent average consensus. Consensus problems are when multiple agents, each with its own dynamics and limited access to the other agents' states, are intended to be stabilized such that all the agents' states are equal. Applications include distributed computing, networks of sensors, flocking, rendezvous, attitude alignment, and synchronization of coupled oscillators; see~\cite{BlWlTc2011},~\cite{WrRwbEma2007} and~\cite{RosRmm2004} and their references. 

Event-triggered strategies for consensus problems have been studied quite extensively over the last decade, with some of the earliest works appearing in 2009~\cite{DvdEf2009,DvdKhj2009,EkXwNh2010}. 
We refer to~\cite{CnEgJc2019} for a detailed survey on the history of this problem but summarize the relevant points next.
 
A seminal work on this topic is~\cite{DvdEfKhj2012}, which develops centralized event- and self-triggered strategies that lead to multi-agent consensus, and then modifies them to be distributed. Unfortunately, although the centralized event-triggered strategy is able to guarantee a positive MIET, the distributed strategies are unable to guarantee the prevention of Zeno solutions. 
Similarly, the results in~\cite{EgYcHyPaDc2013} are unable to exclude Zeno behavior. 

Some works have addressed this issue by considering a periodically sampled (or sampled-data) implementation to trivially address this issue, but this assumes perfect synchronization among the entire network which is neither practical or scalable~\cite{XmTc2013,XmLxYcsCnGjp2015,CnJc2016,AaAaAm2018}. 
More recent works have even considered asynchronous periodic implementations but
these require some sort of global knowledge to find periods that will work~\cite{FxTc2016,YlCnZtQl2017}. We are instead interested in designing event-triggering conditions that guarantee a positive MIET rather than forcing one artificially. 

More related to our work,~\cite{GssDvdKhj2013,BcZl2017,XyKlDvdKhj17} present distributed event-triggered strategies for the consensus problem that prevent Zeno solutions and ensure convergence to consensus. The first two include an explicit function of time in the trigger mechanism, while the third uses a dynamic triggering mechanism, by including a virtual state. 
While these are a good start, unfortunately none of these can guarantee a positive MIET for the agents, which is our main goal. 

The distributed event-triggered strategy in~\cite{VsdMaWpmhh2017} is able to guarantee convergence to consensus with a positive MIET enforced; however, it requires global parameters in order to design each agent's controller, so it is not fully distributed. {Alternatively, and most similarly to the methods used in this work, the authors of~\cite{CdpRp2017} utilize a hybrid systems formulation to solve a closely related problem in which a different communication model is considered.} In their work they show the existence of a positive MIET for a fully distributed event-triggered strategy that guarantees asymptotic convergence with an event trigger that employs a dynamic virtual state; however, this work still requires a type of synchronization as agents need to trigger events in pairs. 


Finally, another important consideration is the robustness of a MIET. {
We note here that we are primarily concerned with the robustness of the event-triggering strategy rather than robustness in terms of feedback stabilization.}
The authors of~\cite{DpbWpmhh2014} acutely point out that, even if an event-triggered controller may guarantee a positive MIET, it is possible that arbitrarily small disturbances remove this property which means it is still equally useless for implementing on physical systems. {We refer to~\cite{DpbWpmhh2014} which introduces a notion of strongly non-Zenoness for a more precise mathematical definition for this type of robustness.} Another potential issue for event-triggered strategies is robustness to imperfect event detection. More specifically, analysis on event-based solutions often rely on the very precise timings of actions in response to events. Consequently, we are also interested in designing a solution that is robust to small timing errors in determining when event conditions have been satisfied.


{
\emph{Statement of contributions:} This problem revisits a simple single-integrator multi-agent average consensus problem originally conceived in~\cite{DvdEf2009,DvdKhj2009,EkXwNh2010}. By re-formulating the problem using hybrid systems, we are able to provide a novel algorithm with several key fundamental improvements over similar event-triggered strategies in the literature. The contributions of this paper are threefold. First, we provide the first known fully distributed solution that guarantees a positive MIET. 
Second, we develop a method to design the triggering functions such that each agent is able to independently prescribe their guaranteed MIET in a distributed way; which has important implications on implementability of the proposed algorithms in real applications. Finally, we investigate the robustness of the MIET against both imperfect event detection, for which we provide an alternative robust trigger, and additive state disturbances, discussing its effects.} Simulations illustrate our results.

\section{Preliminaries}\label{se:preliminaries}

The Euclidean norm of a vector $v \in \real^n$ is denoted by $||v||$. An n-dimensional column vector with every entry equal to $1$ is denoted by $\mathbf{1}_n$, and an n-dimensional column vector with every entry equal to $0$ is denoted by $\mathbf{0}_n$. 
The minimum eigenvalue of a square matrix $A$ is given by $\operatorname{eigmin}(A)$ and its maximum eigenvalue is given by $\operatorname{eigmax}(A)$. 
The distance of $x$ from the set $\mathcal{A}$, which is $\operatorname{min}_a ||x-a||$, where $a \in \mathcal{A}$, is denoted by $||x||_\mathcal{A}$.
Given a vector~$v \in \real^N$, we denote by $\operatorname{diag}(v)$ the $N\times N$ diagonal matrix with the entries of $v$ along its diagonal.

Young's inequality is
\begin{align}\label{eq:young}
xy \leq \dfrac{a}{2}x^2 + \dfrac{1}{2a}y^2,
\end{align}
for $a>0$ and $x,y\in \real$~\cite{GhhJelGp1952}.

By $V^{-1}(C)$ where $C \subset \real^m$, we denote the set of points $\{s \in \real^n : V(s) \in C\}$, for a function $V: \real^n \rightarrow \real^m$.
By $\real_{\geq 0}$ we denote the set of nonnegative real numbers, and by $\integernonnegative$ we denote the set of nonnegative integers.
The closure of a set $U \in \real^n$ is denoted by $\overline{U}$. 
The domain of a mapping $f$ is denoted by $\operatorname{dom}f$ and its range is denoted by $\operatorname{range}f$.

\paragraph*{Graph Theory}
An unweighted graph $\mathcal{G} = (V, \mathcal{E}, A)$ has a set of vertices $V = \{1, 2, ..., N\}$, a set of edges $\mathcal{E} \subset V \times V$, and an adjacency matrix $A \in \real^{N\times N}$ with each entry $a_{ij} \in \{0, 1\}$, where $a_{ij}=1$ if $(i, j) \in \mathcal{E}$, and $a_{ij} = 0$ otherwise. For a digraph (directed graph), edge $(i, j)$ is distinct from edge $(j, i)$.
 A path between vertex $i$ and vertex $j$ is a finite sequence of edges $(i,k)$, $(k,l)$, $(l,m)$, $\dots$, $(n,j)$. A digraph is strongly connected if there exists a path between any two vertices. 
A weighted digraph is one where each edge $(i, j) \in \mathcal{E}$ has a weight $w_{ij} > 0$ to it. For an edge $(i, j)$, $j$ is an out neighbor of $i$ and $i$ is an in neighbor of $j$. The in-degree, $d_i^\text{in}$, for a vertex $i$ is the sum of all the weights for the edges that correspond to its in neighbors, and the out-degree, $d_i^\text{out}$, is the same for its out neighbors. A weight-balanced digraph is a digraph where $d_i^\text{in}=d_i^\text{out} = d_i$ for each vertex $i$. A weighted digraph has a weighted adjacency matrix $A$ where the $ij$th element is the weight for edge $(i,j)$. For a weight-balanced digraph, the degree matrix $D^\text{out}=D^\text{in}$ is a diagonal matrix with $d_i$ as the $i$th diagonal element, and the Laplacian is $L=D^\text{out}-A$.

\paragraph*{Hybrid Systems}
A hybrid system $\mathcal{H} = (C,f,D,G)$ is a tuple composed of a flow set $C\in \real^n$, where the system state $x \in \real^n$ continuously changes according to $\dot{x}=f(x)$, and a jump set $D \in \real^n$, where $x$ discretely jumps to $x^+ \in G(x)$, where $f$ maps $\real^n \rightarrow \real^n$ and $G : \real^n \rightrightarrows \real^n$ is set valued~\cite[Definition 2.2]{RgRgsArt2012}.  While $x \in C$, the system can flow continuously and while $x \in D$, the system can jump discontinuously.

A compact hybrid time domain is a subset $E_\text{compact} \subset \real \times \natural$ for which $E_\text{compact} = \cup_{j=0}^{J-1}([t_j, t_{j+1}],j)$, for a finite sequence of times $0 \leq t_0 \leq t_1 \leq ... \leq t_J$, and a hybrid time domain is a subset $E\subset \real \times \natural$ such that $\forall (T,J)\in E$, $E \cap ([0,T]\times \{0,1,...J\})$ is a compact hybrid time domain~\cite[Definition 2.3]{RgRgsArt2012}. The hybrid time domain is used to keep track of both the elapsed continuous time $t$ and the number of discontinuous jumps $j$. 

See the appendix for more definitions and results relating to hybrid systems.

\section{Problem Formulation}\label{se:statement}

\renewcommand{\arraystretch}{1}
\begin{table}
\begin{center}
\begin{tabular}{|ll|}
\hline
Definition & Domain \\
\hline
$q_i = \left[x_i, \hat{x}_i, \chi_i \right]^T$ & $\in \real^3$ \\ 
$v_i = (q_i, \{ \hat{x}_j \}_{j \in \NN_i^\text{out}})$ & $\in \real^3 \times \real^{|\NN_i^\text{out}|}$ \\
$e_i = x_i - \hat{x}_i$ & $\in \real$ \\
$\hat{z}_i = (L\hat{x})_i = \sum_{j \in \Nouti} w_{ij}(\hat{x}_i - \hat{x}_j)$ & $\in \real$ \\
$\hat{\phi}_i = \sum_{j \in \Nouti} w_{ij}(\hat{x}_i - \hat{x}_j)^2$ & $\in \realnonnegative$ \\
\hline
\end{tabular}
\end{center}
\caption{Agent~$i$ model definitions.}\label{tab:notation}
\end{table}


%
We begin by stating a long-standing version of the event-triggered consensus problem. We then show why existing solutions to it are not pragmatic and how we reformulate the problem to obtain solutions that can be implemented on physical platforms. Please see Table~\ref{tab:notation} for a summary of notation used in the following.

Consider a group of $N$ agents whose communication topology is described by a directed, weight-balanced, and strongly connected graph $\mathcal{G}$ with edges $\mathcal{E}$ and Laplacian matrix $L$. Each agent is able to receive information from its out neighbors and send information to its in neighbors, and each weight of the graph is a gain applied to the information sent from one agent to another. 

The state of each agent~$i$ at time~$t \geq 0$ is given by~$x_i(t)$ with single-integrator dynamics
\begin{align}\label{eq:dynamics}
\dot{x}_i(t) = u_i(t),
\end{align}
where $u_i$ is the input for agent~$i$. It is well known that the input
\begin{align}\label{eq:idealcontrol}
u_i(t) = -\sum_{j \in \Nouti} w_{ij} (x_i(t) - x_j(t) )
\end{align}
drives all agent states to the average of the initial conditions~\cite{RosRmm2004}, which is defined as
\begin{align*}
\bar{x} \triangleq \dfrac{1}{N}\sum_{i=1}^{N}x_i(0).
\end{align*}
Note that under the control law~\eqref{eq:idealcontrol}, the average~$\bar{x}(t)$ is an invariant quantity. Defining $x = [x_1, x_2, \dots, x_N ]^T$ and $u = [u_1,u_2, \dots, u_N]^T$ as the vectors containing all the state and input information about the network of agents, respectively, we can describe all inputs together by
\begin{align*}
u(t) = -Lx(t). 
\end{align*}
However, in order to implement this control law, each agent must have continuous access to the state of each of its out neighbors. Instead, we assume that each agent~$i$ can only measure its own state~$x_i$ and must receive neighboring state information through wireless communication. We consider event-triggered communication and control where each agent only broadcasts its state to its neighbors at discrete instances of time. More formally, letting $\{t_\ell^i \}_{\ell \in \integernonnegative} \subset \realnonnegative$ be the sequence of times at which agent~$i$ broadcasts its state to its in neighbors~$j \in \NN_i^\text{in}$, the agents instead implement the control law
\begin{align}\label{eq:input}
u(t) &= -\hat{z} \triangleq -L\hat{x}(t) ,
\end{align}
where $\hat{x} = [ \hat{x}_1 , \dots ,\hat{x}_N ]^T$ is the vector of the last broadcast state of each agent. More specifically,
given the sequence of broadcast times~$\{t_\ell^i\}_{\ell \in \integernonnegative}$ for agent~$i$, we have 
\begin{align}\label{eq:communication}
\hat{x}_i(t) &= x_i(t^i_\ell) \quad \text{for} \quad t \in [t^i_\ell,t^i_{\ell+1}).
\end{align}
Note that the input~\eqref{eq:input} still ensures that the average of all agent states is an invariant quantity because $\dot{\bar{x}} = \dfrac{1}{N}\mathbf{1}_N^T\dot{x} = \dfrac{1}{N}\mathbf{1}_N^T(-L\hat{x})=0$, which follows from the weight-balanced property of the graph.

At any given time~$t \geq 0$, we define
\begin{align}\label{eq:localmemory}
v_i(t) \triangleq (x_i(t),\hat{x}_i(t),\{\hat{x}_j(t)\}_{j \in \Nouti})
\end{align}
as all the dynamic variables locally available to agent~$i$. The problem of interest, formalized below, is then to obtain a triggering condition based on this information such that the sequence of broadcasting times~$\{t_\ell^i\}_{\ell \in \integernonnegative}$, for each agent~$i$, guarantees that the system eventually reaches the average consensus state. 

\begin{problem}\longthmtitle{Distributed Event-Triggered Consensus}\label{pr:mainoriginal}
Given the directed, weight-balanced, and strongly connected graph $\mathcal{G}$ with dynamics~\eqref{eq:dynamics} and input~\eqref{eq:input}, find a triggering condition for each agent~$i$, which depends only on locally available information $v_i$, such that
$x_i \rightarrow \bar{x}$ for all~$i \in \until{N}$. 
\end{problem}

This problem was first formulated in~\cite{DvdEf2009,DvdKhj2009} in 2009. Since then, there have been many works dedicated to this problem in both the undirected~\cite{XmTc2013,GssDvdKhj2013,DvdEfKhj2012,EgYcHyPaDc2013,XyKlDvdKhj17,JbCn2019} and directed~\cite{XmLxYcsCnGjp2015,CnJc2016,PxCnZt2018} cases. 
Unfortunately, although the above referenced papers provide theoretical solutions to this problem, we are not aware of a single solution which can be implemented on physical systems when considering some practical concerns which are described next. For details on the history of this problem and its many theoretical solutions we refer the interested reader to~\cite{CnEgJc2019}, but we summarize the main points here. 

The earliest solutions to this problem did not adequately investigate the Zeno phenomenon which invalidates their correctness~\cite{DvdEfKhj2012,EgYcHyPaDc2013}. In particular, these solutions to Problem~\ref{pr:mainoriginal} did not rule out the possibility of Zeno behavior meaning that it was possible for a sequence of broadcasting times~$\{t_\ell\}_{\ell \in \integernonnegative}$ to converge to some finite time~$t_{\ell} \rightarrow T > 0$. This is clearly troublesome since all theoretical analysis then falls apart after~$t > T$, invalidating the asymptotic convergence results. 
More recently, the community has acknowledged the importance of ruling out Zeno behavior to guarantee that all the sequences of times~$t_\ell^i \rightarrow \infty$ as~$\ell \rightarrow \infty$ for all~$i \in \until{N}$. While enforcing this additional constraint on the sequences of broadcasting times guarantees that \emph{theoretically} the solutions will converge to the average consensus state, there are still some important practical issues that must be considered.

Even if it can be guaranteed that Zeno behaviors do not occur and the inter-event times are strictly positive for all agents~$i$, 
\begin{align*}
t_{\ell+1}^i - t_{\ell}^i > 0 ,
\end{align*}
unfortunately this is still not enough to guarantee that the solution can be realized by physical devices. This is because although the inter-event times are technically positive, they can become arbitrarily small to the point that no physical hardware exists that can keep up with the speed of actions required by the event-triggered algorithm. The solutions in~\cite{GssDvdKhj2013},~\cite{BcZl2017},~\cite{XyKlDvdKhj17}, and~\cite{PxCnZt2018} have this problem. 
This is inherently different from guaranteeing a strictly positive MIET~$\tau$, where $t_{\ell+1}^i - t_\ell^i \geq \tau > 0$, which is the focus of our work here. 

Specifically, we consider the case where each agent~$i$ has some maximum rate $\frac{1}{\tau_i}$ at which it can take actions (e.g., broadcasting information, computing control inputs). That is, each agent~$i$ cannot broadcast twice in succession in less than~$\tau_i$ seconds. In other words, each agent~$i \in \until{N}$ is limited by hardware in terms of how fast they are able to take actions,
\begin{align*}
t_{\ell+1}^i - t_\ell^i \geq \tau_i,
\end{align*}
for all~$\ell \in \integernonnegative$.
 Note that there are also solutions that guarantee a MIET, but make other sacrifices to do so. For example, the solution in~\cite{GssDvdKhj2013} is able to guarantee a MIET under certain conditions, but convergence is only to a neighborhood of consensus. 
Additionally, the algorithms in~\cite{VsdMaWpmhh2017,AaAaAm2018} are able to enforce a MIET, but only by using global parameters of the system to design the algorithm, which is impractical in cases where the parameters may change or are otherwise difficult to measure. 
The algorithm in~\cite{CdpRp2017} is fully distributed and has a positive MIET, but it still requires pair of agents to trigger events at the same time, which necessitates synchronization. 
With all this in mind, we reformulate Problem~\ref{pr:mainoriginal} such that solutions to the problem can be implemented on physical platforms given that each agent~$i$ is capable of processing actions at a frequency of up to~$\frac{1}{\tau_i}$.


\begin{problem}\longthmtitle{Distributed Event Triggered Consensus with Designable MIET}\label{pr:main}
Given the directed, weight-balanced, and strongly connected graph $\mathcal{G}$ with dynamics~\eqref{eq:dynamics}, input~\eqref{eq:input}, and the minimum periods~$(\tau_1, \dots, \tau_N)$ for each agent, find a triggering
condition for each agent $i$, which depends only on local information $v_i$, such that
$x_i \rightarrow \bar{x}$ and 
\begin{align}\label{eq:MIETcondi}
\min_{\ell \in \integernonnegative} t^i_{\ell + 1} - t^i_\ell \geq \tau_i ,
\end{align}
for all~$i \in \until{N}$.
\end{problem}

{To the best of our knowledge, Problem~\ref{pr:main} has not yet been fully solved. Rather than being able to guarantee 
a strictly positive minimum inter-event time, similar works often settle for only ruling out Zeno executions meaning communication between agents may still need to occur arbitrarily fast in order for the convergence results to hold. Conversely, other works more simply force a MIET through the use of a dwell-time at the cost of losing the exact asymptotic convergence guarantee in exchange for practical consensus. The works~\cite{GssDvdKhj2013,BcZl2017,XyKlDvdKhj17} only preclude Zeno behavior but cannot guarantee that~\eqref{eq:MIETcondi} holds even for arbitrarily small minimum periods~$\tau_i$. The algorithm proposed in~\cite{VsdMaWpmhh2017} comes close but requires global system information. The methodology used in~\cite{CdpRp2017} is similar to ours here but ultimately solves a different problem. Thus, we provide the first complete solution to Problem~\ref{pr:main}.} For now we consider no state disturbances and perfect event detection but will relax these assumptions in Section~\ref{se:robustness} where we study the robustness of the {event-triggering strategy} with respect to various forms of disturbances.

\subsection{Hybrid Systems Formulation}
In order to solve Problem~\ref{pr:main}, we first reformulate it using hybrid systems tools; similar to~\cite{CdpRp2017}. {We refer to the original state~$x$ as the `physical' state that represents the actual state that we wish to control. Separately, we maintain a set of `cyber' or virtual states corresponding to the internal memory of each agent. Given the communication model described by~\eqref{eq:communication}, it seems natural to keep track of the last broadcast state~$\hat{x}_i$ for each agent as one of the cyber states. Additionally, we introduce an extra virtual state $\clovar_i$ for each agent~$i\in\until{N}$ to introduce dynamics into our triggering strategy and collect these components in the vector~$\clovar = [\clovar_1, \clovar_2, \dots, \clovar_N]^T$. Note that the internal variable~$\clovar_i$ is only available to agent~$i$. In this work we consider scalar internal variables~$\clovar_i \in \real$ but note that more sophisticated controllers or even learning-based controllers could be captures by increasing the complexity of the internal variables. The hybrid systems formulation will aid us here in properly modeling both the continuous-time dynamical system with discrete-time memory and control updates.}

{We now define the extended state vector for agent~$i$ as
\begin{align*}
q_i(t) = \left[ \begin{array}{c}
x_i(t)\\
\hat{x}_i(t)\\
\clovar_i(t)
\end{array} \right],
\end{align*} 
and the extended state (capturing both `physical' and `cyber' states) of the entire system is $q = [q_1^T, q_2^T, \dots, q_N^T]^T \in \real^{3N}$. 

With a slight abuse of notation, we now redefine the local information~\eqref{eq:localmemory} available to agent~$i$ at any given time as its own extended state~$q_i$ and the last broadcast states of its out-neighbors~$\{\hat{x}_j\}_{j \in \NN_i^\text{out}}$,
%
\begin{align}\label{eq:localmemory2}
v_i(t) \triangleq (q_i(t), \{ \hat{x}_j(t) \}_{j \in \NN_i^\text{out}}) . 
\end{align}
\
}

In this work, the goal is to use the internal variable~$\clovar_i$ to determine exactly when agent~$i$ should broadcast its current state to its neighbors. To achieve this we let the state take values~$\clovar_i \geq 0$ and prescribe an event-trigger whenever~$\clovar_i = 0$. The exact dynamics of~$\clovar_i$ will be designed later in Section~\ref{se:design} to guarantee a solution to Problem~\ref{pr:main}.
%
More specifically, at any given time~$t \geq t_\ell^i$, the next triggering time~$t_{\ell+1}^i$ is given by
\begin{align}\label{eq:trigFunc}
t^i_{\ell+1} = \inf \{t \geq t^i_{\ell}: \clovar_i(t) = 0 \text{ and }\hat{x}_i \neq x_i \}, 
\end{align}
for all~$\ell \in \integernonnegative$, for each agent $i$. 

%
%
With the role of the new internal variable~$\clovar_i$ established (although its dynamics will be designed later), we can now formalize our hybrid system
\begin{align}\label{eq:system} 
\mathcal{H} = (C,f,D,G).
\end{align}
We refer to the appendix and~\cite{RgRgsArt2012} for formal hybrid systems concepts and assume the reader is familiar with this formalism.

\noindent \textbf{Flow set} (between event-triggering times):
The flow set~$C$ for the entire system is then given by
\begin{align}\label{eq:flowset}
C = \{q \in \real^{3N}: \clovar_i \geq 0 \text{ for all } i \in \until{N} \}.
\end{align}
While the system state~$q \in C$, the system flows according to~$f$
\begin{align}\label{eq:flowdynamics}
\dot{q} = f(q) = \left[ \begin{array}{c}
f_1(v_1)\\
\vdots\\
f_N(v_N)\\
\end{array} \right] \quad \text{for } q \in C,
\end{align}
with the individual extended states evolving according to
\begin{align}
\dot{q}_i = f_i(v_i) \triangleq \left[ \begin{array}{c}
-\hat{z}_i\\
0\\
\gamma_i(v_i)
\end{array} \right],
\end{align}
where $(L\hat{x})_i$ denotes the $i$th element of the column vector $L\hat{x}$ and~$\gamma_i$ is the function to be designed. 

The first row is exactly~$\dot{x}_i~=~u_i$ as defined in~\eqref{eq:input}, the second row says the last broadcast state is not changing between event-triggers~$\dot{\hat{x}_i}~=~0$, and the last row is the dynamics of the internal variable~$\dot{\clovar_i}~=~\gamma_i(v_i)$. Recall that~$v_i$ defined in~\eqref{eq:localmemory2} only contains information available to agent~$i$ and thus the dynamics of the internal variable~$\clovar_i~=~\gamma_i(v_i)$ must be a function only of~$v_i$.


\noindent \textbf{Jump set} (at event-triggering times):
The jump set~$D$ for the entire system is then given by
\begin{align}\label{eq:jumpset}
D = \cup_{i=1}^N \{q \in \real^{3N}: \clovar_i \leq 0 \}.
\end{align}
Although in a general hybrid system there may be no notion of distributed information, since in this work the jumps correspond to \emph{some} agent triggering an event, we formally define
\begin{align}\label{eq:jumpsetlocal}
D_i \triangleq \{q \in \real^{3N}: \clovar_i \leq 0 \}
\end{align}
as the subset of~$D$ corresponding to when specifically agent~$i$ is responsible for the jump. 
{For~$q \in D_i$, we consider the following local jump map
\begin{align*}
g_i(q) = \left[ \begin{array}{c} q_1^+ \\ \vdots \\ q_i^+ \\ \vdots \\ q_{N}^+ \end{array} \right] \triangleq \left[ \begin{array}{c} q_1 \\ \vdots \\ \left( \begin{array}{c} x_i \\ x_i \\ \clovar_i\end{array} \right) \\ \vdots \\ q_{N} \end{array} \right].
\end{align*}
More specifically, letting~$t_\ell^i$ be the time at which agent~$i$ triggers its~$\ell$th event~$q(t_\ell^i) \in D_i$, this map leaves the physical state and the dynamic variable unchanged~$x_i^+ = x_i,\chi_i^+ = \chi_i$, and updates its ``last broadcast state'' to its current state~$\hat{x}_i^+ = x_i$.} Note also that this leaves all other agents' states unchanged~$q_j^+ = q_j$ for all~$j \neq i$.


Now, since multiple agents may trigger events at once, the jump map must be described by a set-valued map~$G : \real^{3N} \rightrightarrows \real^{3N}$~\cite{RgsJjbbNvdwWpmhh2014,CdpRp2017}, where
\begin{align}\label{eq:jumpDyn}
G(q) \in \{\dots, g_i(q), \dots\},
\end{align}
for all $i$ such that $q \in D_i$. Note that this construction of the jump map ensures that it is outer-semicontinuous, which is a requirement for some hybrid systems results.


Now, we reformulate Problem~\ref{pr:main} in a more structured manner by using the hybrid system~\eqref{eq:system}. More specifically, by using this formulation we have formalized the objective of finding a local triggering strategy to designing the function~$\gamma_i$ that depends only on the local information~$v_i$ defined in~\eqref{eq:localmemory2}.

\begin{problem}\longthmtitle{Distributed Event-Triggered Consensus with Designable MIET}\label{pr:reform}
Given the directed, weight-balanced, and strongly connected graph $\mathcal{G}$ with dynamics~\eqref{eq:dynamics}, input~\eqref{eq:input}, and the minimum periods~$(\tau_1, \dots, \tau_N)$ for each agent, find the dynamics of the clock-like variable, $\gamma_i(v_i)$, such that $x_i \rightarrow \bar{x}$ and 
\begin{align*}
\min_{\ell \in \integernonnegative} t^i_{\ell + 1} - t^i_\ell \geq \tau_i ,
\end{align*}
for all~$i \in \until{N}$.
\end{problem}

\section{Dynamic Event-Triggered Algorithm Design}\label{se:design}
In order to solve Problem~\ref{pr:reform}, we perform a Lyapunov analysis to design~$\gamma \triangleq \dot{\clovar}$, the dynamics of $\clovar$.
Inspired by~\cite{CdpRp2017}, we use a Lyapunov function with two components: $V_P$ represents the physical aspects of the system, while $V_C$ represents the cyber aspects, related to communication and error. 
For convenience, let $e \triangleq x - \hat{x}$ denote the vector containing the error for each agent's state, which is the difference between the actual state and the last broadcast state. We begin by considering 
{
\begin{align*}
V_{P}(q) &= (x - \bar{x})^T(x - \bar{x}) = ||x - \bar{x}||^2\\
V_C(q) &= \sum_{i=1}^N \clovar_i.
\end{align*}}
Note $V_C \geq 0$ because $\clovar_i \geq 0$ for all $i \in \until{N}$. We then consider the Lyapunov function
\begin{align}\label{eq:LyapFunc}
V(q) = V_{P}(q) + V_C(q) .
\end{align}
Note that $V(q) \geq 0$ is continuously differentiable for all $q \in \real^{3N}$. {Moreover,
~$V(q)=0$ when all agents have reached their target state and each clock-like variable~$\clovar_i$ is equal to 0.} 

Now we will examine the evolution of $V$ along the trajectories of our algorithm to see under what conditions it is nonincreasing, and design~$\gamma$ accordingly. 
{In order to do so, we will have to split $\dot{V}$ into components $\dot{V}_i$ such that $\dot{V} = \sum_{i=1}^N \dot{V}_i$ and each $\dot{V}_i$ depends only on the local information~$v_i$ available to agent~$i$. Choosing $V$ properly to ensure that $\dot{V}$ can be split like this is essential to designing $\gamma_i$ and doing so is nontrivial.}
{
Recalling our system flow dynamics~\eqref{eq:flowdynamics}, we write
\begin{align*}
\dot{V}_P &= -2(x - \bar{x})^T\hat{z}\\
\dot{V}_C &= \sum_{i=1}^N \gamma_i.
\end{align*}
Because the graph is weight-balanced, $\bar{x}^TL = \mathbf{0}_N^T$. Therefore, $\dot{V}_{P} = -2x^T\hat{z} = -2\hat{x}^T\hat{z} -2 e^T\hat{z}$ and
\begin{align*}
\dot{V} = \dot{V}_{P} + \dot{V}_C = -2\hat{x}^T\hat{z} - 2e^T\hat{z} + \sum_{i=1}^N \gamma_i.
\end{align*}
Expanding this out and using notation defined in Table~\ref{tab:notation} yields
\begin{align}\label{eq:V1dot}
\dot{V} &= \sum_{i=1}^N \left( -\sum_{j \in \Nouti} w_{ij}(\hat{x}_i - \hat{x}_j)^2 - 2e_i(\hat{z}_i  + \gamma_{i} \right) \\
\dot{V} &= \sum_{i=1}^N \dot{V}_i = \sum_{i=1}^N \left( -\hat{\sumvar}_i - 2e_i\hat{z}_i + \gamma_{i} \right) \nonumber.
\end{align}}

{
We are now interested in designing~$\gamma_i$ for each agent~$i \in \until{N}$ such that~$\dot{V}(q) \leq 0$. 
Therefore, we choose
\begin{align}\label{eq:clockdyn1}
\gamma_i = \sigma_i\hat{\phi}_i + 2e_i\hat{z}_i,
\end{align}
where $\sigma_i \in (0,1)$ is a design parameter. Note that whenever an agent~$i$ triggers an event the error~$e_i$ is immediately set to 0, and since~$\dot{\clovar}_i = \gamma_i \geq 0$ at these times we ensure that~$\clovar_i \geq 0$ at all times~$t \geq 0$. We can now write the derivative of the Lyapunov function as
\begin{align}\label{eq:vdot}
\dot{V} = \sum_{i=1}^N -(1-\sigma_i)\hat{\phi}_i \leq 0.
\end{align}}
{This choice of the clock-like dynamics~$\gamma_{i}$ is continuous in $q$ for constant $\hat{x}$ }{and ensures that~$\dot{V} \leq 0$.}

\subsection*{Algorithm Synthesis}
With the dynamics~$\gamma_i$ of the clock-like variable~$\clovar_i$ defined for each agent~$i \in \until{n}$, we can now summarize all the components of our synthesized distributed dynamic event-triggered coordination algorithm and formally describe it from the viewpoint of a single agent.

The control input at any given time~$t \geq 0$ is 
\begin{align*}
u_i(t) = -\hat{z}(t)_i = -\sum_{j \in \Nouti}w_{ij} (\hat{x}_i(t)-\hat{x}_j(t)).
\end{align*}
{The sequence of event times~$\{t_\ell^i\}_{\ell \in \integernonnegative}$ at which agent~$i$ broadcasts its state to neighbors is given by each time
the clock-like variable reaches zero when that agent's error is nonzero, i.e., 
\begin{align*}
t_\ell^i = \inf \{t \geq t_{\ell - 1}^i : \clovar_i(t) = 0 \text{ and } e_i \neq 0\}.
\end{align*}
}
The algorithm is formally presented in Table~\ref{tab:algorithm}.

\begin{table}[htb]
  \centering
  \framebox[.9\linewidth]{\parbox{.85\linewidth}{%
    \parbox{\linewidth}{Initialization; at time~$t = 0$ each agent~$i \in \until{N}$ performs:}
      \vspace*{-2.5ex}
      \begin{algorithmic}[1]
      \STATE Initialize~$\hat{x}_i = x_i$
      \STATE Initialize~$\clovar_i = 0$
      \end{algorithmic}
      \parbox{\linewidth}{At all times $t$ each agent $i \in \until{N}$
        performs:}
      \vspace*{-2.5ex}
      \begin{algorithmic}[1]
        \IF{
          $\clovar_i = 0$ and $e_i \neq 0$}  \label{algo:1}
        \STATE set $~\hat{x}_i = x_i$ (broadcast state information to neighbors)
        \STATE set $u_i = -\sum_{j \in \Nouti} w_{ij}(\hat{x}_i - \hat{x}_j)$ (update control signal)
        \ELSE
        \STATE propagate~$\clovar_i$ according to its dynamics~$\gamma_i$ in~\eqref{eq:clockdyn1}
        \ENDIF
        \IF{
          new information $\hat{x}_k$ is received from some neighbor(s)~$k \in \Nouti$} 
        \STATE update control signal~$u_i = -\sum_{j \in \Nouti} w_{ij}(\hat{x}_i - \hat{x}_j)$
        \ENDIF
      \end{algorithmic}}}
  \caption{Distributed Dynamic Event-Triggered Coordination Algorithm.}\label{tab:algorithm}
\end{table}

\section{Main Results}
Here we present the main results of the paper by discussing the properties of our algorithm. We begin by finding the guaranteed positive minimum inter-event time (MIET) for each agent and showing how it can be tuned individually. 

\begin{theorem}[Positive MIET]\label{th:ldtResult}
Given the hybrid system $\mathcal{H}$, if each agent~$i$ implements the distributed dynamic event-triggered coordination algorithm presented in Table~\ref{tab:algorithm} {with~$\sigma_i \in (0,1)$,
 then the inter-event times for agent~$i$ are lower-bounded by
\begin{align}\label{eq:MIET1}
T_i \triangleq \frac{\sigma_i}{d_i} > 0.
\end{align}
}
That is,
\begin{align*}
t_{\ell+1}^i - t_\ell^i \geq T_i
\end{align*}
for all~$i \in \until{N}$ and~$\ell \in \integernonnegative$.
\end{theorem}
\begin{IEEEproof}
See the appendix.
\end{IEEEproof}

\begin{remark}[Design Trade off]\label{re:tradeoff}
{\rm
The design parameter $\sigma_i$ represents the trade off between larger inter-event times and faster convergence speeds. Larger $\sigma_i$ makes the magnitude of $\dot{V}$ smaller~\eqref{eq:vdot}. However, it is a coefficient of the nonnegative term in $\chi_i$'s dynamics~\eqref{eq:clockdyn1}, so larger $\sigma_i$ means longer inter-event times, and increasing it increases the MIET~\eqref{eq:MIET1}. 
Additionally, note that the MIETs can be guaranteed up to a maximum of {$T_{i,\text{max}} \triangleq \dfrac{1}{d_i}$ for each agent $i$, because $\sigma_i \in (0,1)$.} 
}
\end{remark}

Next, we present our main convergence result. To the best of our knowledge, this is the first work to design a fully distributed event-triggered communication and control algorithm that guarantees asymptotic convergence to the average consensus state with a lower bound on the agent-specific MIET that can be chosen by the designer.


\begin{theorem}[Global Asymptotic Convergence]\label{th:mainResult}
Given the hybrid system~$\mathcal{H}$, if each agent~$i$ implements the distributed dynamic event-triggered coordination algorithm presented in Table~\ref{tab:algorithm} with agent $i$ triggering events when 
$\clovar_i = 0$ {and $\hat{x}_i \neq x_i$ with~$\sigma_i \in (0,1)$, then all trajectories of the system are guaranteed to asymptotically converge to the set
\begin{align*}
\{q : \hat{\phi}_i = 0\text{ } \forall\text{ } i\}.
\end{align*}}
\end{theorem}
\begin{IEEEproof}
See the appendix.
\end{IEEEproof}

{
\begin{remark}[Convergence]\label{rm:conv}
{\rm
From Theorem~\ref{th:mainResult}, the algorithm presented in Table~\ref{tab:algorithm} does not entirely solve Problem~\ref{pr:reform}, because agents only converge to $\{q : \hat{\phi}_i = 0\text{ } \forall\text{ } i\}$, and, in order for $x_i \rightarrow \bar{x}$ for each agent~$i$, we also require  and $e_i = 0$, $\forall i$. 
Therefore, if the algorithm is modified with an additional trigger which guarantees that each agent $i$ will always broadcast again, eventually, when $e_i \neq 0$, then full convergence is guaranteed. 
More formally, for $\ell \in \mathbb{Z}_{\geq 0}$, if we guarantee that $t^i_{\ell+1}<\infty$ when $\exists t \in [t^i_\ell,t^i_{\ell+1})$ such that $e(t) \neq 0$, then $x_i(t) \rightarrow \bar{x}$ as~$t \rightarrow \infty$, for all~$i \in \until{N}$. 

A very simple way to do so is to set a maximum time $T_\text{max}^i \geq T_i$ between events for each agent $i$, so that, if it triggers an event at time $t_\ell^i$, the latest another event can trigger is $t_\ell^i + T_i$.
}
\end{remark}
}

\subsection*{Minimum Inter-Event Time Design}
{As noted in Remark~\ref{re:tradeoff}, the design parameter $\sigma_i$ can be used to choose the desired MIET for agent $i$, up to a maximum of $T_{i,\text{max}}$.} According to Problem~\ref{pr:main}, we must be able to guarantee that the lower-bound on the MIET~$T_i$ as provided in Theorem~\ref{th:ldtResult} is greater than or equal to the prescribed~$\tau_i$. 

Consequently, if~$\tau_i < T_{i,\text{max}}$ for all~$i \in \until{N}$, then it is easy to see how~\eqref{eq:MIET1} in Theorem~\ref{th:ldtResult} can directly be used to choose the design {parameter} appropriately for each agent. 
{{In the case that there exists some agent(s)~$j$ such that~$\tau_j \geq T_{j,\text{max}}=\frac{1}{d_j}$,
 then the graph must be redesigned so that $d_j$ is lower for each agent $j$. Note that there exist distributed methods of choosing these gains for an existing strongly connected digraph so that it will be weight-balanced~\cite{BgJc2012,AirTcCnh2014}. However, further analysis here is beyond the scope of this paper.}

\section{Robustness}\label{se:robustness}
{A problem with many event triggered algorithms is a lack of robustness
guarantees in the triggering strategies. In particular, we consider the effects of two different types of disturbances that are often problematic for event-triggered control systems, discussing how our algorithm can be implemented to preserve the MIET in the presence of additive state disturbances and providing a modification that makes it robust against imperfect event detection.} 

\subsection*{Robustness Against State Disturbances}
We first analyze the robustness of our MIET against state disturbances. As noted in~\cite{DpbWpmhh2014}, simply guaranteeing a positive MIET may not be practical if the existence of arbitrarily small disturbances can remove this property, resulting again in solutions that might require the agents to take actions faster than physically possible in an attempt to still ensure convergence. Therefore, it is desirable for our algorithm to exhibit robust global event-separation as defined in~\cite{DpbWpmhh2014}, which means that the algorithm can still guarantee a positive MIET for all initial conditions even in the presence of state disturbances, which is referred to as a robust MIET.

{
More formally, we desire the MIET given in~\eqref{eq:MIET1} to hold, even in the presences of arbitrary disturbances. Instead of the deterministic dynamics~\eqref{eq:dynamics}, consider
\begin{align}\label{eq:distributedynamics}
\dot{x}_i(t) = u_i(t) + w_i(t),
\end{align}
where~$w_i(t)$ is an arbitrary, unknown, additive state disturbance applied to each agent's state. However, ensuring that the MIET holds in these circumstances depends on the specific implementation, as discussed in the following remark.

\begin{remark}[Trigger Robustness]\label{rm:robTrig}
{\rm
Note that, following an event at time $t_\ell^i$, each agent $i$ does not even need to check the trigger condition~\eqref{eq:trigFunc} until time $t_\ell^i+T_i$, because, by Theorem~\ref{th:ldtResult}, it cannot be satisfied before that time, in the absence of any disturbance. Therefore, in implementation, each agent can wait $T_i$ seconds after triggering an event before triggering a new one. This will have no effect on the performance of the algorithm in the absence of disturbances, but it will ensure that the MIET is observed in the presence of disturbances of the form given in~\eqref{eq:distributedynamics}. 
Additionally, see the definition of the hybrid system $\mathcal{H}'$ in~\eqref{eq:modSys} for how this can be modeled in theory.
} \oprocend
\end{remark}
Note that this does not guarantee convergence all the way to consensus in the presence of disturbances, simply that the positive MIET will be preserved.} 
Analyzing the actual convergence properties in any formal sense is beyond the scope of this work. 
Instead, we consider the following simple example to show how the algorithm may handle a disturbance. 
If each $w_i$ is an independent and identical Gaussian process with zero mean and variance $\sigma^2$, then dynamics of the average position, $\dot{\bar{x}}$, will be a random variable with $E[\dot{\bar{x}}] = 0$ and $\text{var}(\dot{\bar{x}}) = \sigma^2/N$. This indicates that $\bar{x}$ is a Wiener process, which has a Gaussian distribution with a mean equal to the initial average and a variance of $t\frac{\sigma^2}{N}$.  {We demonstrate the effects of such a disturbance in Section~\ref{se:simulations}.}

\subsection*{Robustness Against Imperfect Event Detection}
In addition to robustness against state disturbances, another important source of uncertainty that cannot be overlooked in event-triggered control systems
is imperfect event detection. Event-triggered controllers are generally designed and analyzed assuming very precise timing of different actions is possible while continuously monitoring the event conditions. This is not only impractical but problematic if not accounted for in the event-triggered control design. 

The algorithm presented in Table~\ref{tab:algorithm} is no longer guaranteed to converge if there are delays in the triggering times. More specifically, let~$t_{\ell+1}^{i^*}$ be the time at which the event condition~\eqref{eq:trigFunc} is actually satisfied, but the condition is not actually detected until the actual triggering time~$t_{\ell+1}^i = t_{\ell+1}^{i^*} + \delta t^i_\ell$, where~$\delta t^i_\ell \in [0,\Delta_i]$ is a random delay in the detection of the event with an upper-bound of~$\Delta_i > 0$. We refer to this as triggering the $\ell$th event, for agent~$i$, $\delta t^i_\ell$ seconds late. An example where this would be the case is if the trigger condition for agent~$i$ is checked periodically, with frequency $1/\Delta_i$.

Fortunately, our algorithm is still guaranteed to converge as long as agent~$i$ triggers before~$\chi_i$ reaches~$0$. By utilizing the knowledge of the maximum delay~$\Delta_i$,  we can modify the local jump set~\eqref{eq:jumpsetlocal} that defines when agent~$i$ should trigger an event such that we are guaranteeing that all events are detected by the time~$\chi_i$ reaches 0 when accounting for the delay. Let
\begin{align*}
\tilde{D}_i = \{q \in \mathbb{R}^{3N} : h_i(v_i) \in [0, \Delta_i]\},
\end{align*}
where
\begin{align*}
h_i(v_i) = \left\lbrace
\begin{matrix}
\frac{\sigma_i}{d_i}\left( 1 - \frac{e_i^2}{\chi_i+e_i^{2}}\right), & \text{for } (\chi_i,e_i) \neq (0,0)\\
\frac{\sigma_i}{d_i}, & \text{otherwise}
\end{matrix}
  \right. .
\end{align*}

As will be shown in the proof, this provides a time window with a length of $\Delta_i$ seconds starting from when $h_i(v_i) = \Delta_i$ during which agent~$i$ can trigger an event and still guarantee $\chi_i \geq 0$. Therefore, we define a robust event trigger condition
\begin{align}\label{eq:modTrigger}
h_i(v_i) \leq \Delta_i .
\end{align}
Note that the MIET with this new event-trigger condition is now $\frac{\sigma_i}{d_i} - \Delta_i$. 
}
This is formalized next.

{
\begin{theorem}[Robust Convergence with MIET]\label{th:robMIET}
Given the hybrid system $\mathcal{H}$, if each agent~$i$ implements the distributed dynamic event-triggered coordination algorithm presented in Table~\ref{tab:algorithm} and each event is triggered at most $\Delta_i$ seconds after~\eqref{eq:modTrigger} is satisfied, with $\sigma \in (0,1)$ and $\Delta_i \in [0, \frac{\sigma_i}{d_i})$, then all trajectories of the system are guaranteed to asymptotically converge to the set
\begin{align*}
\{q : \hat{\phi}_i = 0\text{ } \forall\text{ } i\},
\end{align*}
for all~$i \in \until{N}$. 
Additionally, the inter-event times for agent~$i$ are lower-bounded by
\begin{align}\label{eq:robMIET1}
\tilde{T}_{i} \triangleq & \frac{\sigma_i}{d_i} - \Delta_i > 0.
\end{align}
\end{theorem}}
\begin{IEEEproof}
See the appendix.
\end{IEEEproof}

\begin{remark}[Trigger Robustness]\label{rm:robust}
{\rm
The implication of Theorem~\ref{th:robMIET} is the following. 
Intuitively, rather than agent~$i$ waiting for the condition~\eqref{eq:trigFunc} to be satisfied exactly and respond immediately, it simply begins triggering an event {when condition~\eqref{eq:trigFunc} could be satisfied soon, and as long as the event can be detected and fully responded to before then, the algorithm will work as intended.} However, note that this imposes a trade off because triggering earlier will result in a shorter guaranteed MIET, as shown in the result of Theorem~\ref{th:robMIET}.
} \oprocend
\end{remark}

\section{Simulations}\label{se:simulations}
\begin{figure*}
\subfigure[]{\includegraphics[width=.3\linewidth]{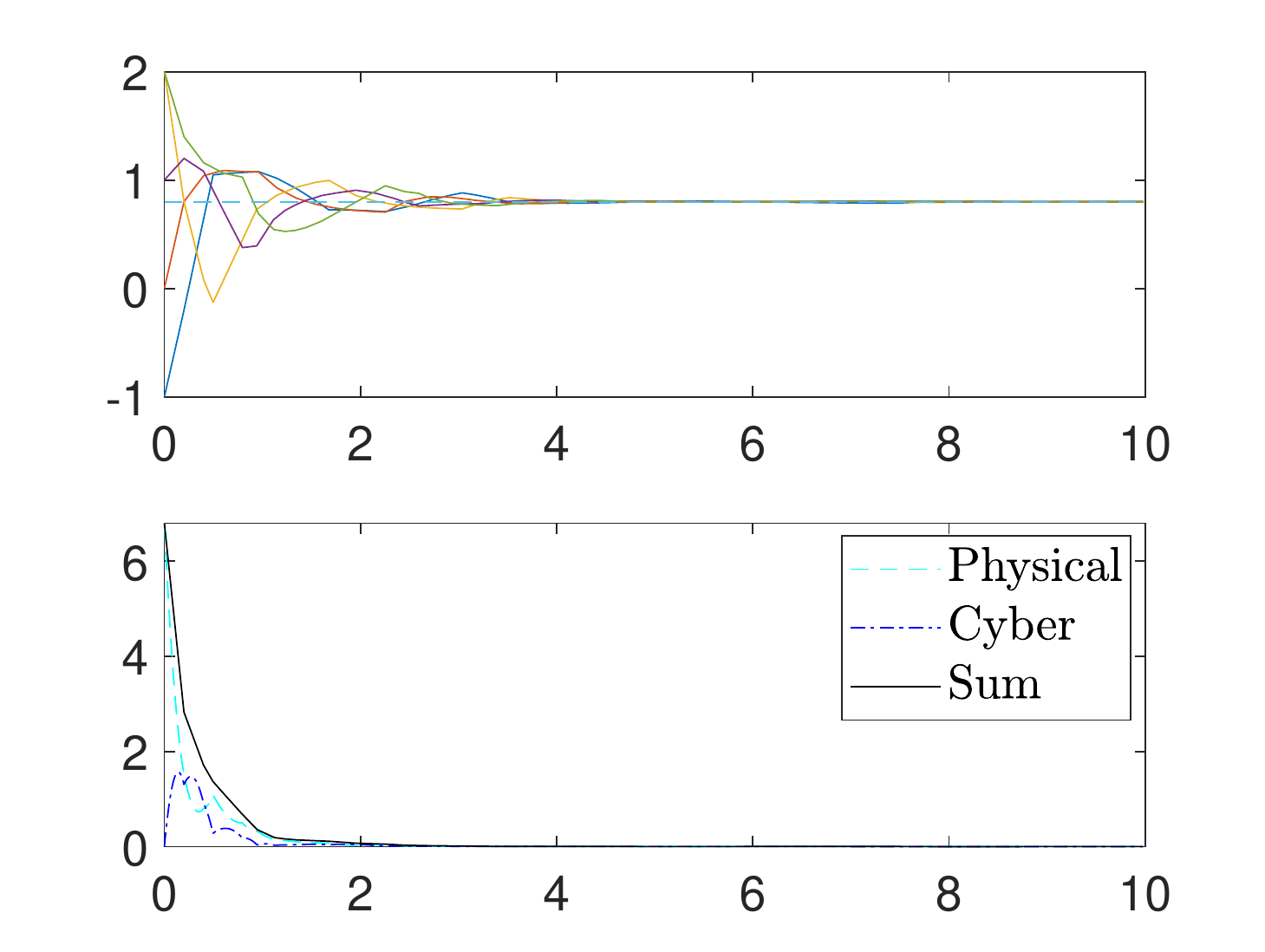}
\put(-77,-3){{\small $t$}}%
\put(-150,87){{\small $x$}}%
\put(-150,29){{\small $V$}}
} \hfill
\subfigure[]{\includegraphics[width=.3\linewidth]{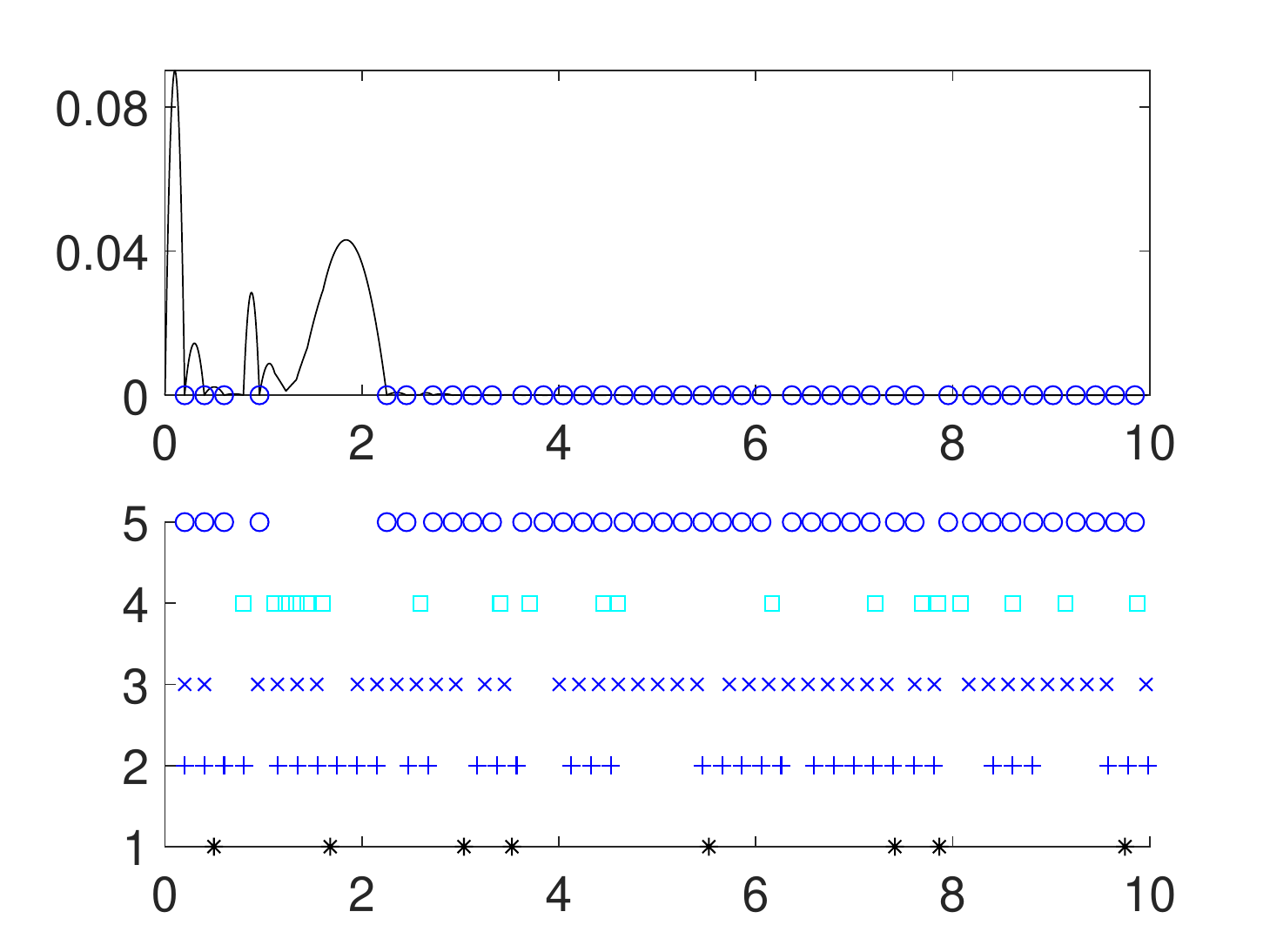}
\put(-77,-3){{\small $t$}}%
\put(-160,87){{\small $\chi_5$}}%
\put(-168,29){{\small Events}}
} \hfill
\subfigure[]{\includegraphics[width=.3\linewidth]{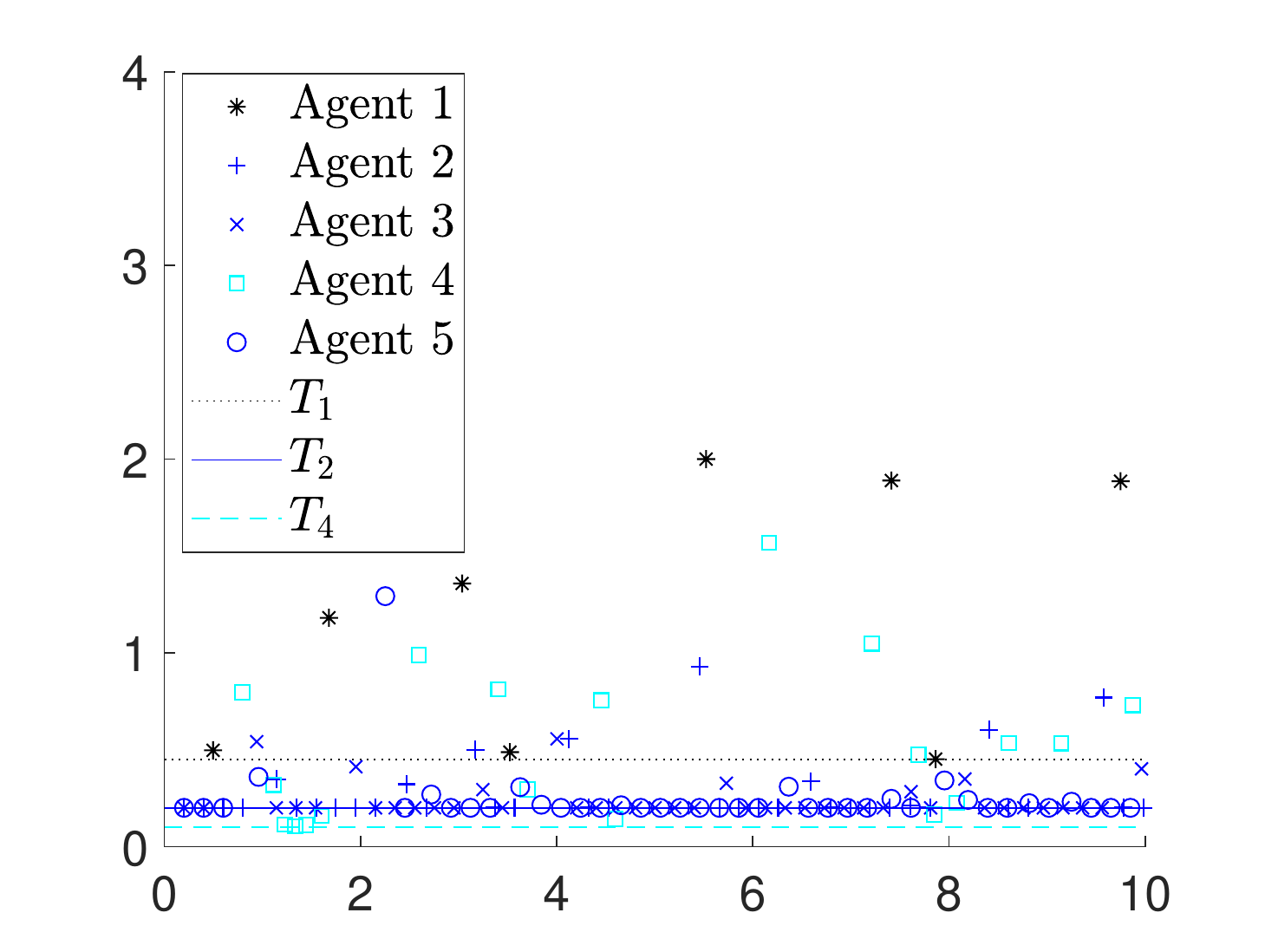}
\put(-77,-3){{\small $t$}}%
\put(-175,61){{\tiny $t^i_{\ell+1}-t_\ell^i$}}%
} \hfill
\caption{{Plots of the simulation results of the dynamic event-triggered algorithm showing (a) the trajectories of the agents (top), with the dashed line representing the average, and the evolution of the whole Lyapunov function $V$ (bottom) as well as the physical component $V_P$ and the cyber component $V_C$; (b) the clock-like state variable~$\clovar_5$ for agent~$5$ (top) and rows of stars indicating the event times of all agents (bottom); and (c) the inter-event times~$t^i_{\ell+1}-t_\ell^i$ of each agent~$i$, with the lower bounds as computed by~\eqref{eq:MIET1} marked by the lines.}
}\label{fig:TrajectGraphs}
\end{figure*}

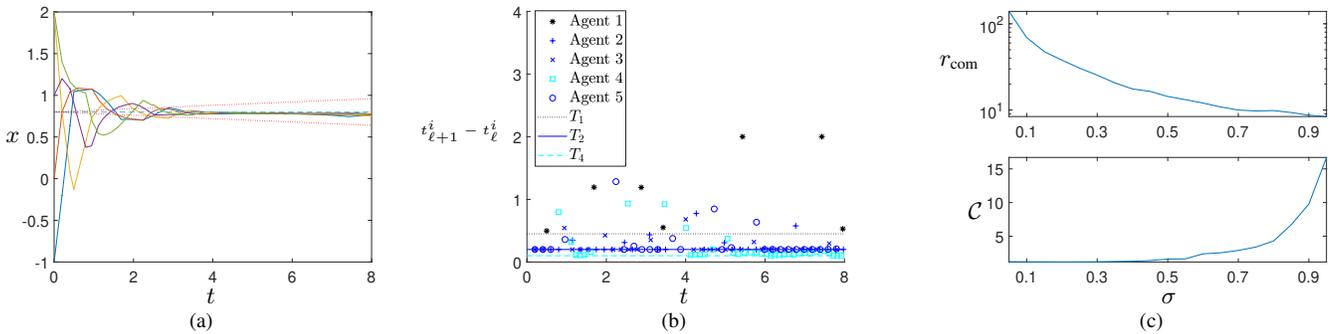
\begin{figure*}
\subfigure[]{\includegraphics[width=.3\linewidth]{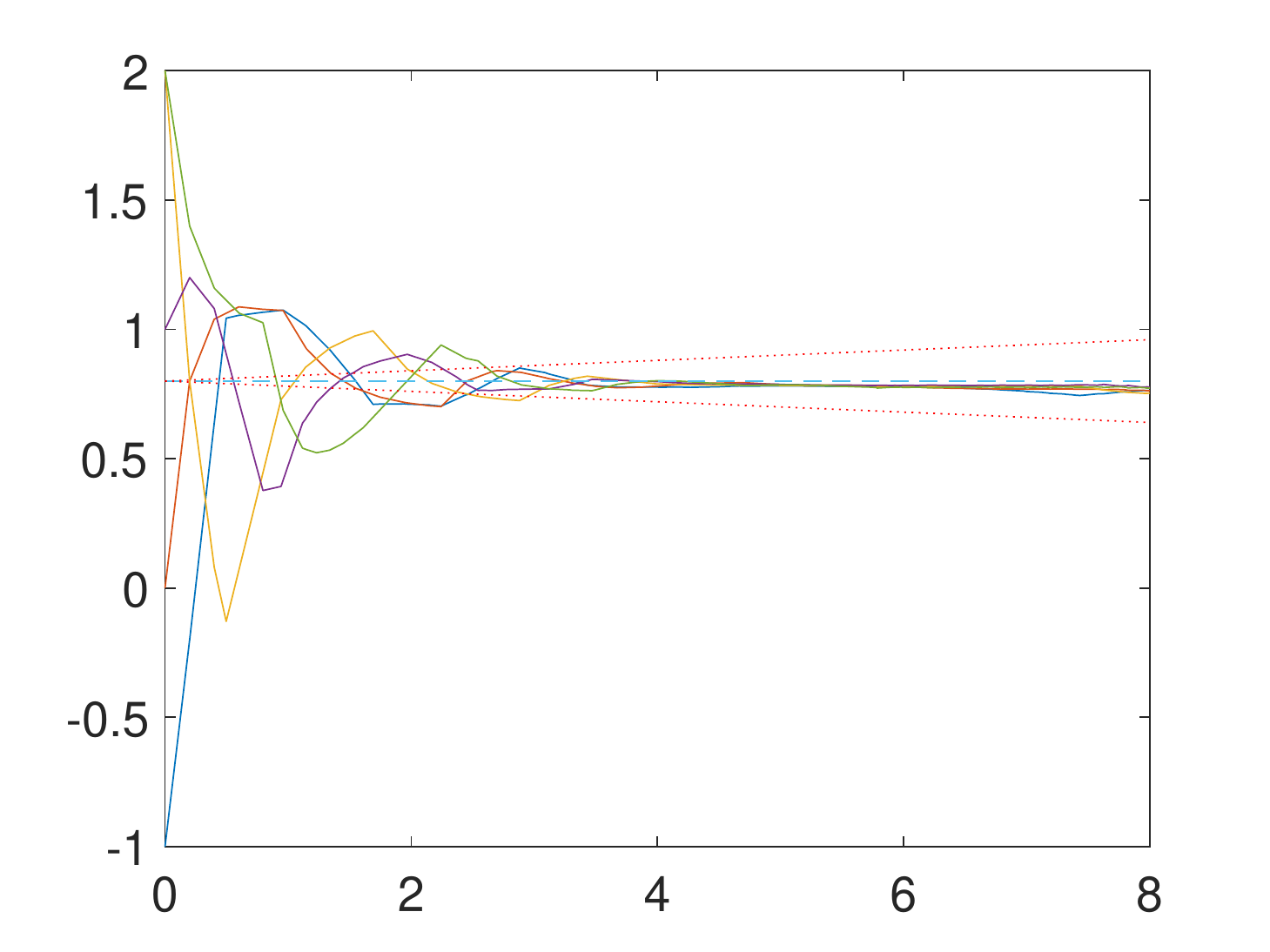}
\put(-77,-3){{\small $t$}}%
\put(-153,58){{\small $x$}}%
} \hfill
\subfigure[]{\includegraphics[width=.3\linewidth]{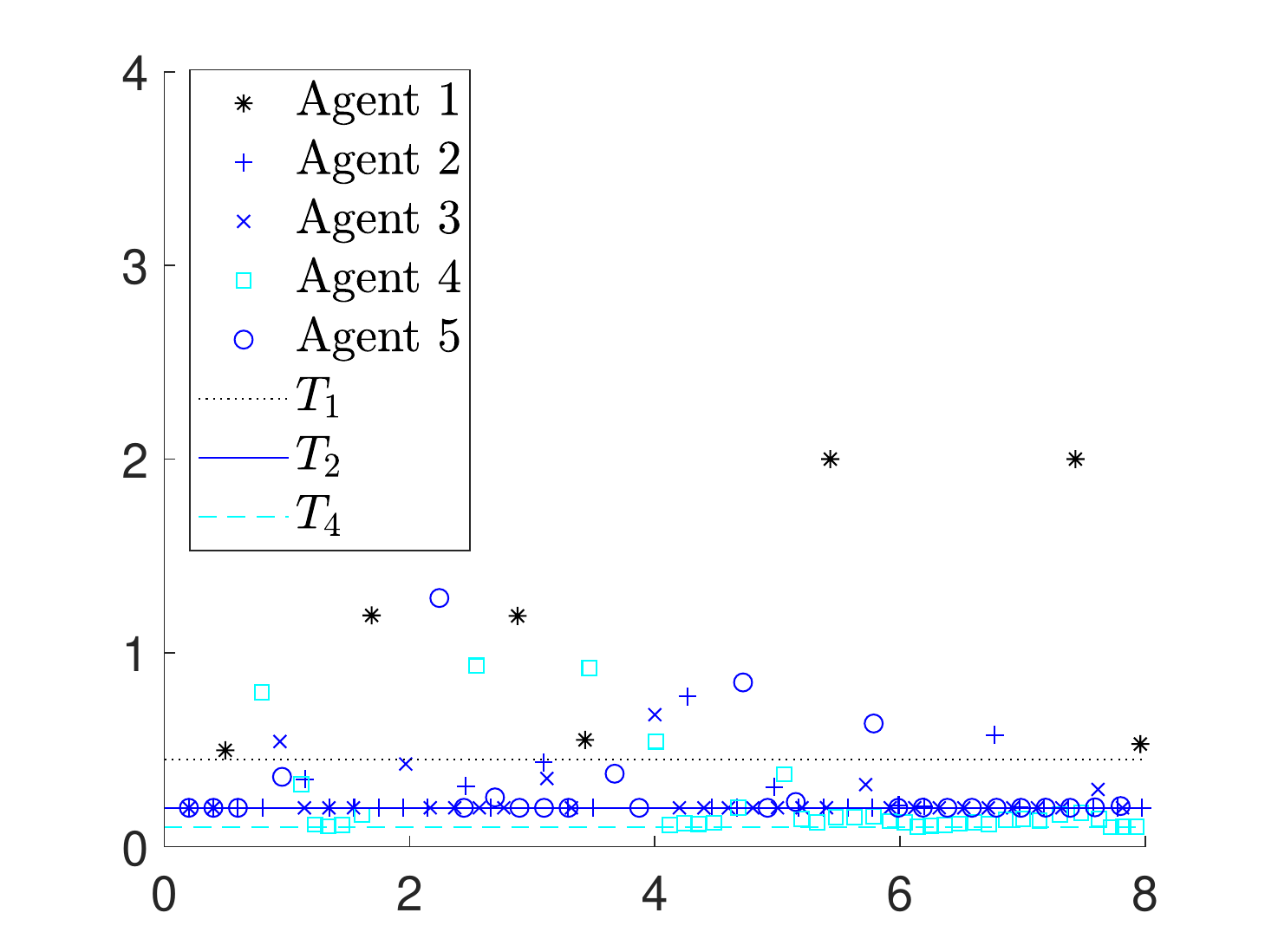}
\put(-77,-3){{\small $t$}}%
\put(-175,61){{\tiny $t^i_{\ell+1}-t_\ell^i$}}%
} \hfill
\subfigure[]{
\includegraphics[width=0.3\linewidth]{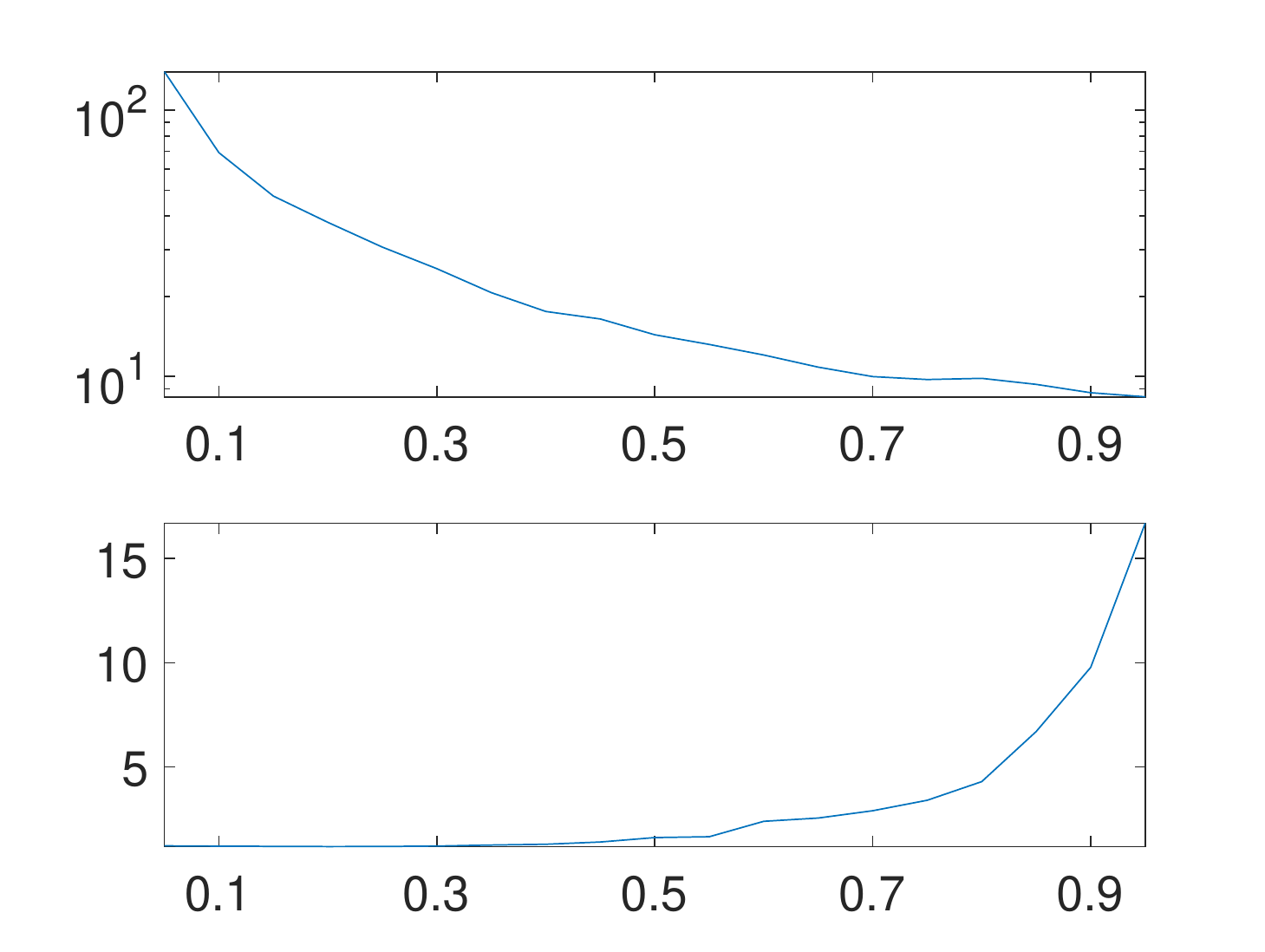}
\put(-77,-3){{\small $\sigma$}}%
\put(-161,87){{\small $r_\text{com}$}}%
\put(-150,29){{\small $\mathcal{C}$}}
}%
\caption{{ Simulation results with additive state disturbances showing (a) the trajectories of the agents subjected to zero-mean additive white Gaussian noise with a variance of $0.1$ (the dashed blue line indicates the expected value of the average position ($\bar{x}(0)$), and the dotted red lines show the variance over time ($\bar{x}(0)\pm t\frac{\sigma_w^2}{N}$); (b) the inter-event times~$t^i_{\ell+1}-t_\ell^i$ of each agent~$i$, with the lower bounds as computed by~\eqref{eq:MIET1} marked by the lines; and (c) the communication rate (top) and~$\mathcal{H}_2$-norm cost~\eqref{eq:cost} (bottom) as we vary the design parameter~$\sigma$.}}\label{fig:varyParms}
\end{figure*}

{
To demonstrate our distributed event-triggered control strategy, we perform various simulations using~$N = 5$ agents and a directed graph whose Laplacian is given by 
\begin{align*}
L = \left[ \begin{array}{ccccc}
 2 & -1 &  0 &  0 & -1\\
 0 &  2 &  0 &  0 & -2\\
-2 &  0 &  2 &  0 &  0\\
 0 & -1 & -2 &  3 &  0\\
 0 &  0 &  0 &  -3 & 3
\end{array} \right].
\end{align*}

Additionally, we consider that agents have a minimum operating period of $[\tau_1,\tau_2,\tau_3,\tau_4,\tau_5]=[0.4,0.25,0.25,0.1,0.2]$. Because the agents have out degrees of $[d_1,d_2,d_3,d_4,d_5]=[2,2,2,3,3]$, we use~\eqref{eq:MIET1} and choose $[\sigma_1,\sigma_2,\sigma_3,\sigma_4,\sigma_5] = [0.9,0.4,0.4,0.3,0.6]$, so that the guaranteed MIETs for the agents are $[T_1,T_2,T_3,T_4,T_5] = [0.45,0.2,0.2,0.1,0.2]$.

All simulations use the same initial conditions of $\hat{x} = x = [-1, 0, 2, 1, 2]^T$ in order to explore the effects of the different design parameters. 

The top plot in Figure \ref{fig:TrajectGraphs} shows the main results. Figure \ref{fig:TrajectGraphs} (a) shows the positions of the agents over time, demonstrating that they converge to initial average position, indicated by the dashed line, and the bottom plot shows the evolution of the Lyapunov function, which can be seen to be nonincreasing, although the physical portion and the cyber portion are allowed to increase individually. 
The top plot in Figure \ref{fig:TrajectGraphs} (b) shows the evolution of the clock-like state, $\clovar_5$, for agent $5$, and the bottom plot shows when each agent triggers an event, demonstrating the asynchronous, aperiodic nature of event triggering.  
Figure~\ref{fig:TrajectGraphs} (c) shows the inter-event times for all agents, with the horizontal lines indicating the theoretical lower bounds. For each agent $i$, the minimum inter-event time $T_i$ can be seen to be respected, and the bound appears to be tight, as expected from the theoretical analysis.

Figure \ref{fig:varyParms} (a) shows the effect of applying an additive white Gaussian noise disturbance, i.e. $\dot{x} = u + w$, where each element of $w$ is an independent and identically distributed Gaussian process, with zero mean and a variance of $\sigma_w^2 = 0.1$. To ensure that the MIET is respected in the presence of noise, the algorithm is implemented in a self-triggered fashion. That is, instead of measuring $e_i$ to propagate the dynamics~\eqref{eq:clockdyn1}, $e_i$ is approximated assuming no noise. This simulation suggests that the expected value of each agent's state is the current average position, although that average can now change with time. Figure \ref{fig:varyParms} (b) shows that the minimum inter-event times are still respected with this implementation.

Next, to show the effect of the design parameter $\sigma_i$ on the algorithm's performance, we set $\sigma_i = \sigma$ for each agent~$i$ and varied it.
Figure \ref{fig:varyParms} (c) shows the results on 2 statistics: the average communication rate $r_\text{com}$, which is the number of events divided by the simulation length, and the cost $\mathcal{C}$, defined as follows. Considering each agent's difference from the average as an output, similar to~\cite{SdYgNm2018}, we adopt the square of the $\mathcal{H}_2$-norm of the system as a cost performance metric 
\begin{align}\label{eq:cost}
\mathcal{C} \triangleq \int_{t=0}^{t=T_\text{max}}\sum_{i=1}^N \left(x(t) - \bar{x}\right)^2,
\end{align}
where $T_\text{max} = \infty$. However, for simplicity in simulations, we use the rough approximation of $T_\text{max}=20$ being the simulation length. The choice of parameter $\sigma$ can be seen to be a trade off between cost (speed of convergence) and communication rate. Higher values of $\sigma$ result in a higher cost, but also requires less communication and results in a higher MIET by~\eqref{eq:MIET1}.}



\section{Conclusions}\label{se:conclusions}
This paper has used the multi-agent average consensus problem to present a dynamic agent-focused event-triggered mechanism which ensures stabilization and prevents Zeno solutions by allowing for a chosen minimum inter-event time for each agent. The algorithm is fully distributed in that it not only requires no global parameters, but the correctness of the algorithm can also be guaranteed by each agent individually. That is, no global conditions (besides connectivity of the graph) need to even be checked to ensure the overall system asymptotically converges. Additionally, it provides robustness against missed event times, guaranteeing convergence as long as events are triggered within a certain window of time. 

While this work has presented an algorithm that distributed agents can implement to guarantee asymptotic convergence, further research is needed to study the transient properties or our proposed and related algorithms. {We plan to examine this algorithm to see if it can guarantee exponential convergence to consensus, and, in particular, how the secondary trigger discussed in Remark~\ref{rm:conv} should be designed for good performance. }


\appendix
{The necessary results from hybrid systems for the rigorous proof are presented first, then a modified hybrid system which can account for different event triggers is presented, and then a general convergence result is given in Lemma~\ref{lem:genConverge}. Finally, the proofs of Theorems~\ref{th:ldtResult},~\ref{th:mainResult}, Lemma~\ref{lem:timeToEvent}, and~\ref{th:robMIET} follow.}

\subsection{Hybrid Systems Results}
\begin{lemma}[Lemma 5.10 in~\cite{RgRgsArt2012}]\label{th:OSCcondi}
A set-valued mapping $M : \real^n \rightrightarrows \real^m$ is outer semicontinuous if and only if the graph of $M$ is closed.
\end{lemma}
\begin{definition}[Definition 2.4 in~\cite{RgRgsArt2012}]
A function $\solvar$: $E \rightarrow \real^n$ is a hybrid arc if $E$ is a hybrid time domain and if for each $j \in \natural$, the function $t \rightarrow \solvar(t,j)$ is locally absolutely continuous on the interval $I^j = \{t: (t,j) \in E\}$.
\end{definition}

\begin{definition}[Definition 2.6 in~\cite{RgRgsArt2012}]
A hybrid arc $\solvar$ is a solution to the hybrid system $(C,F,D,G)$ if $\solvar \in \overline{C} \cup D$, and
\begin{enumerate}
\item for all $j \in \natural$ such that $I^j \triangleq \{t: (t,j) \in \operatorname{dom}\solvar\}$ has a nonempty interior
\begin{align*}
\solvar(t,j) \in C\ \forall\ t \in \operatorname{int} I^j
\end{align*}
\begin{align*}
\dot{\solvar}(t,j) = F\left(\solvar(t,j)\right) \text{ for almost all } t \in I^j
\end{align*}
\item for all $(t,j) \in \operatorname{dom} \solvar$ such that $(t,j+1)\in \operatorname{dom}\solvar$
\begin{align*}
\solvar(t,j) \in D
\end{align*}
\begin{align*}
\solvar(t,j+1) \in G\left(\solvar(t,j)\right).
\end{align*}
\end{enumerate}
\end{definition}

\begin{definition}[Definition 2.7 in~\cite{RgRgsArt2012}]
A solution $\solvar$ is maximal if there does not exist another solution $\solvar '$ to $\mathcal{H}$ such that $\operatorname{dom}\solvar$ is a strict subset of $\operatorname{dom} \solvar '$ and $\solvar(t,j) = \solvar '(t,j)\ \forall\ (t,j)\in \operatorname{dom}\solvar$.
\end{definition}

\begin{theorem}[Theorem 6.8 in~\cite{RgRgsArt2012}]\label{th:nomWellPosed}
If a hybrid system $\mathcal{H}$ satisfies the following assumption, then it is nominally well-posed.
$C$ and $D$ are closed subsets of $\real^n$;

\begin{enumerate}
\item $F:\real^n \rightarrow \real^n$ is outer semicontinuous and locally bounded relative to $C$, $C$ is a subset of the domain of $F$, and $F(q)$ is convex for every $q \in C$;
\item $G:\real^n \rightrightarrows \real^n$ is outer semicontinuous and locally bounded relative to $D$ and $D$ is a subset of the domain of $G$.
\end{enumerate}
\end{theorem}
A solution $\solvar$ is complete if $\text{dom}\ \solvar$ is unbounded.
\begin{definition}[Weak Invariance~\cite{RgRgsArt2012}]\label{def:weakInvariance}
{\rm
Given a hybrid system $\mathcal{H}$, a set $S \subset \real^n$ is said to be

\begin{itemize}
\item \emph{weakly forward invariant} if, for every $q \in S$, $\exists$ a complete solution $\solvar$ to $\mathcal{H}$ with initial condition $q$ whose range is a subset of $S$,

\item \emph{weakly backward invariant} if for every $q \in S$ and every $\tau > 0$, there exists at least one maximal solution $\solvar$ to $\mathcal{H}$ with initial condition in $S$ such that for some $(t^*,j^*)$ in the domain of $\solvar$, $t^* + j^*>\tau$, it is the case that $\solvar(t^*,j^*)=q$ and $\solvar(t,j) \in S\ \forall\ (t,j)$ in the domain of $\solvar$ with $t+j\leq t^* + j^*$,

\item and \emph{weakly invariant} if it is both weakly forward invariant and weakly backward invariant.
\end{itemize}
}
\end{definition}
For a solution $\solvar$ to a hybrid system $\mathcal{H}$, $t(j)$ denotes the least time $t$ such that $(t,j)$ is in its domain and $j(t)$ denotes the least index $j$ such that $(t,j)$ is in its domain~\cite{RgRgsArt2012}.
Given $V:\real^n \rightarrow \real$, any functions $u_C, u_D : \real^n \rightarrow [-\infty, \infty]$, and a set $U \subset \real^n$, it is said that the growth of $V$ along solutions to $\mathcal{H}$ is bounded by $u_C, u_D$ on $U$ if for any solution $\solvar$ to $\mathcal{H}$ with its range in $U$,
\begin{align}\label{eq:boundedGrowth}
V\left(\solvar(\overline{t}, \overline{j})\right) - V\left(\solvar(\underline{t}, \underline{j})\right) & \leq \int_{\underline{t}}^{\overline{t}} u_C\left(\solvar \left(s,j(s)\right)\right)ds \nonumber \\
& + \sum_{j=\underline{j}+1}^{\overline{j}}u_D\left(\solvar \left(t(j),j\right)\right)
\end{align}
for all $(\underline{t}, \underline{j}), (\overline{t}, \overline{j})$ in the domain of $\solvar$ such that $(\underline{t}, \underline{j}) \prec (\overline{t}, \overline{j})$~\cite{RgRgsArt2012}.
Intuitively, $u_C$ acts as a bound for the growth of $V$ during system flow, while $u_D$ acts as a bound for its growth during system jumps.

\begin{theorem}\longthmtitle{Invariance Principle for Hybrid Systems~\cite{RgRgsArt2012}}\label{th:invariance}
Consider a continuous function $V :\real^n \rightarrow \real$, any functions $u_C, u_D : \real^n \rightarrow [-\infty, \infty]$, and a set $U \subset \real^n$ such that $u_C(q), u_D(q) \leq 0$ $\forall\ q \in U$ and such that the growth of $V$ along solutions to $\mathcal{H}$ is bounded by $u_C, u_D$ on $U$. Let a complete, bounded solution to $\mathcal{H}$, $\solvar^*$, be such that the closure of its range $\in U$. Then, for some $r \in V(U)$, $\solvar^*$ approaches the nonempty set that is the largest weakly invariant subset of
\begin{align}\label{eq:invarianceSet}
\mathcal{S} \triangleq V^{-1}(r) \cap U \cap \left[\overline{u_C^{-1}(0)} \cup \left(u_D^{-1}(0) \cap G\left(u_D^{-1}(0)\right)\right)\right].
\end{align}
\end{theorem}

\subsection{Proofs}

{
\emph{More General Hybrid System Formulation}

To aid in several proofs, we define another hybrid system which will account for the different event-triggers discussed in this paper. This is a more generalized version of~\eqref{eq:system}, with the jump set extended, using a new timer state, to account for earlier events. We begin by extending the state with a timer variable $\mathcal{T}_i$ for each agent~$i$, so that, with a slight abuse of notation, we have
\begin{align*}
q_i = \left[ \begin{array}{c}
x_i\\
\hat{x}_i\\
\clovar_i\\
\mathcal{T}_i
\end{array} \right],
\end{align*}
$q = [q_1^T,\hdots,q_n^T]^T$, and $\dot{q} = f'(q)$, where $f'(q) = [f_1^{'T},\hdots,f_n^{'T}]^T$ and
\begin{align*}
\dot{q}_i = f_i'(q) \triangleq \left[ \begin{array}{c}
-(L\hat{x})_i\\
0\\
\gamma_i(v_i)\\
1
\end{array} \right] \quad \text{for } q \in C',
\end{align*}
where $C' = \cap _{i=1}^N C'_i$ and
\begin{align*}
C'_i = \{q \in \real^{4N}: \clovar_i \geq 0 \}.
\end{align*}
The extended jump set is now given by~$D' = \cup_{i=1}^N D'_i$, where
{
\begin{align}\label{eq:jumpset}
D'_i =  \{q \in \real^{4N}: \mathcal{T}_i \geq T_i \},
\end{align}}
for any $T_i \in (0,\frac{\sigma_i}{d_i}]$. The new local jump map, for $q \in D_i'$, is
\begin{align*}
g_i'(q) = \left[ \begin{array}{c} q_1^+ \\ \vdots \\ q_i^+ \\ \vdots \\ q_{N}^+ \end{array} \right] \triangleq \left[ \begin{array}{c} q_1 \\ \vdots \\ \left( \begin{array}{c} x_i \\ x_i \\ \clovar_i \\ 0\end{array} \right) \\ \vdots \\ q_{N} \end{array} \right].
\end{align*}
The jump map is described by a set-valued map~$G' : \real^{4N} \rightrightarrows \real^{4N}$
, where
\begin{align}\label{eq:modJumpDyn}
G'(q) \in \{ g_1(q), \dots, g_N(q) \}.
\end{align}
Informally, this defines the timer  $\mathcal{T}_i$ to count up between agent $i$'s events, and it is set to $0$ at each of agent $i$'s events. This means that $\mathcal{T}_i$ measures the time since agent $i$'s most recent event. The extended jump set now allows events to occur as soon as an inter-event time $T_i$ (which can be shorter than the time given in~\eqref{eq:MIET1}, but not longer) has passed, but the system does not have to jump until it is about to exit the flow set when $\exists$ $i:\chi_i = 0$. Finally, we define
\begin{align}\label{eq:modSys}
\mathcal{H}' = (C',f',D',G').
\end{align}

We use this general system~\eqref{eq:modSys} to give the following lemma, which will be used to prove the main results.
\begin{lemma}[General Asymptotic Convergence]\label{lem:genConverge}
Given the hybrid system~$\mathcal{H}'$, for any $q(0) \in \mathbb{R}^{4N}$ such that $x(0) \in \mathbb{R}^{N}$, $\hat{x}(0) = x(0)$, $\chi(0) = \mathbf{0}_N$, and $\mathcal{T} = \mathbf{0}_N$, the system state is guaranteed to asymptotically converge to the set
\begin{align*}
\mathcal{B} \triangleq \{q \in \mathbb{R}^{4N} : \hat{\phi}_i = 0\text{ } \forall\text{ } i\}.
\end{align*}
\end{lemma}}
\begin{IEEEproof}

{Ultimately we wish to apply Theorem~\ref{th:invariance} to our hybrid system. 

Intuitively, we are interested in showing that~$V$ is nonincreasing along the trajectories of our system so that we can apply an invariance principle to show convergence. 
While the state is flowing~($q \in C$), we have already shown that the clock defined by~\eqref{eq:clockdyn1} ensures that~$\dot{V} < 0$ $\forall q \notin \mathcal{B}$. 
When the system jumps~($q \in D$), we have
\begin{align*}
 V \left( g(q) \right) - V(q) = 0,
\end{align*}
for $q \in D$, because $V$ does not depend on $\hat{x}$ or $\mathcal{T}$.

This means that while the system is flowing but not in the target state~($q \in C \setminus \mathcal{B}$), we have~$\dot{V} < 0$. When the state jumps~($q \in D$), the value of~$V$ remains unchanged.} Combining this with the fact that $\bar{x}$ is constant and with the positive MIET result 
to ensure that~$t \rightarrow \infty$ without exhibiting Zeno behavior guarantees {that~$q(t) \rightarrow \mathcal{B}$.} {This is the intuitive argument, which we will formalize next.} 

{
More formally, we are interested in showing that~$V$ is nonincreasing along the trajectories of our system, then characterizing the largest invariant subset where it is zero as~$\mathcal{B}$.

We first show that~$\mathcal{H'}$ is nominally well-posed using~\cite[Theorem~6.8]{RgRgsArt2012}. 
The flow set $C$ and the jump set $D$ are both closed subsets of $\real^{4N}$, $f'(q)$ and $G'(q)$ are outer semicontinuous, and both are locally bounded so all conditions of~\cite[Theorem~6.8]{RgRgsArt2012} are satisfied and the hybrid system $\mathcal{H'}$ is nominally well-posed. 

Now we establish bounds for $V$'s rate of change. 
Let
\begin{align}\label{eq:ud}
u_D = \begin{cases}
 0 & \text{for } q\in D\\
 -\infty & \text{otherwise}\\
\end{cases}.
\end{align}
and
\begin{align}\label{eq:uc}
u_C = \begin{cases}
\sum_{i=1}^N -(1 -\sigma_i)\hat{\phi}_i  & \text{for } q\in C\\
 -\infty & \text{otherwise}\\
\end{cases}.
\end{align}
The function $u_D$ acts as an upper bound on the rate of change of $V(q)$ for each jump (because $V$ does not jump) while $u_C$ acts as an upper bound on the rate of change during flow. 
Note that we have shown that the growth of $V(q)$ along any solution is bounded by $u_C$, $u_D$.

Let $U(\bar{x}(0))\subset \real^{4N} \triangleq \{q \in \real^{4N} : \dfrac{1}{N} \sum_{i=1}^N x_i = \bar{x}(0)\}$. Note that, because the average position of the agents remains constant along all solutions, $U$ is invariant and any solution $\solvar$ such that $\solvar(0) \in U(\bar{x}(0))$ remains in $U(\bar{x}(0))$ for as long as it is defined.

Therefore, by Theorem~\ref{th:invariance}, every complete, bounded solution $\solvar$ such that $\solvar(0) \in U(\bar{x}(0))$ approaches the largest weakly invariant subset of
\begin{align}\label{eq:invarianceSet}
\mathcal{S} \triangleq V^{-1}(r) \cap \left[\overline{u_C^{-1}(0)} \cup \left(u_D^{-1}(0) \cap G\left(u_D^{-1}(0)\right)\right)\right],
\end{align}
for some $r \in V(U(\bar{x}(0)))$. 
The weakly invariant subset of $\overline{u_C^{-1}(0)}$ is $\mathcal{B}$.

As an aside, note that, if we have an addition trigger that guarantees that agent~$i$ will trigger an event and cause a jump in finite time if $e_i \neq 0$, then this set is instead $\{q : \hat{\phi}_i = e_i = 0\text{ } \forall\text{ } i\}$. 

The points in $G(u_D^{-1}(0)) \cap u_D^{-1}(0)$ are points inside the jump set which can be reached by jumping. Note that the timer variables $\mathcal{T}_i$ preclude degenerate cases where the same agent broadcasts more than once at a single time instant, and so the system must leave the jump set after a finite number of jumps. 
Therefore, the largest weakly invariant subset of $\mathcal{S}$ cannot include any points outside $\mathcal{B}$, and so all complete, bounded solutions to $\mathcal{H}'$ starting in $U(\bar{x}(0))$ converge to $\mathcal{B}$ by Theorem~\ref{th:invariance}.

Because $V$ is radially unbounded with respect to $x$ and $\chi$ and $\dot{V} \leq 0$, every solution must be bounded, and, because $C' \cup D' = \real^{4N}$, every solution is complete. Therefore every solution to $\mathcal{H}'$ such that $x(0) \in \real^N$, $\hat{x}(0) = x(0)$, and $\clovar_i \geq 0 $, for $i = 1, 2, \dots , N$ is complete and bounded and converges to $\mathcal{B}$.
}
\end{IEEEproof}

\emph{Proof of Theorem~\ref{th:ldtResult}}

{
To determine the minimum inter-event time for agent $i$, we write the relevant states as their own local dynamical system, and, considering unknowns as inputs to this local system, we apply optimal control to see how quickly the system can be driven to the next event state. The relevant states are the auxiliary variable $\chi_i$ 
and the error $e_i$. 
Since there are no guarantees on the last broadcast states of our neighbors, $\hat{x}_j$ for $j \in \mathcal{N}^\text{out}_i$, we consider the directed distance to each out neighbor as an input to this system. We define our local nonlinear system as
\begin{align}\label{eq:stateDyn}
\dot{\zeta} = f(\zeta,\mu) \triangleq \left[ \begin{matrix}
\sigma_i \mu^T\mathcal{W}\mu + 2\zeta_2W\mu\\
-W\mu\\
\end{matrix}\right],
\end{align}
where $\zeta_1 = \chi_i$, $\zeta_2 = e_i$, $\mu$ is a column vector such that $\mu_j = \hat{x}_i - \hat{x}_j$, $W$ is a row vector such that $W_j = w_{ij}$, and $\mathcal{W} = \text{diag }(W)$. 

To determine the MIET, we want to see how short the time between events can be. Therefore, we look to minimum time optimal control. Our performance index is
\begin{align}\label{eq:perfIndex}
J = \int^{T_i}_01dt,
\end{align}
where $T_i$ is the minimum time we wish to find. 
We assume that the agent has just triggered an event and communicated, so the initial condition is
\begin{align}\label{eq:initCondi}
\zeta(0) = \left[ \begin{matrix}
0\\
0
\end{matrix}\right],
\end{align}
because $\chi_i = \zeta_1 = 0$ is required to trigger the broadcast which sets $e_i = \zeta_2 = 0$. We wish the terminal condition to be when  the condition for a triggering an event is first reached again, that is, when $\zeta_1 = 0$ and $\zeta_2 \neq 0$. However, this set is open, so we cannot use it to perform optimization. Instead, we assume that the final value of $\zeta_2$ is known a priori, so that $\zeta_2(T_i) = e_T$. We can then investigate the effects of different choices of $e_T$. Therefore, the terminal condition is
\begin{align}\label{eq:termCondi}
\zeta(T_i)= \left[ \begin{matrix}
0\\
e_T
\end{matrix}\right].
\end{align}
The Hamiltonian is
\begin{align*}
H &= 1 + \lambda^Tf(\zeta,\mu)\\
&= 1 + \lambda_1(\sigma \mu^T\mathcal{W}\mu + 2\zeta_2W\mu) - \lambda_2W\mu,
\end{align*}
where $\lambda \in \real^2$ is the costate. 
The costate equation~(\cite{FllDlvVls2012}) is
\begin{align*}
-\dot{\lambda} &= \frac{\partial H}{\partial \zeta}^T = \left[\begin{matrix}
0\\
2\lambda_1W\mu
\end{matrix}\right].
\end{align*}
The stationarity condition~\cite{FllDlvVls2012} is
\begin{align*}
0 &= \frac{\partial H}{\partial \mu} = (2\sigma\mathcal{W}\mu + 2\zeta_2W^T)\lambda_1 - W^T\lambda_2.
\end{align*}

Therefore, the relevant solutions for the state, costate, and input can be found to be
\begin{align}\label{eq:optimalInput}
\zeta_2(t) &= -\frac{W\mathcal{W}^{-1}W^T\lambda_2(0)}{2\lambda_1\sigma}t\nonumber\\
\lambda_1 &= \lambda_1(0)\nonumber\\
\lambda_2(t) &= \lambda_2(0) - \frac{W\mathcal{W}^{-1}W^T\lambda_2(0)}{\sigma}t\nonumber\\
\mu &= -\mathbf{1}_{|\mathcal{N}_i|}\frac{e_T}{d_i T_i},
\end{align}
where $\mathcal{W}$ is invertible because it is a diagonal matrix with positive entries on its diagonal. Note that $\mathcal{W}^{-1}W^T = \mathbf{1}_{|\mathcal{N}_i|}$ and so $W\mathcal{W}^{-1}W^T = \sum_{j\in \mathcal{N}_i}w_{ij} = d_i$. Note that this $\mu$ will be useful in other proofs.

However, $T_i$ is still unknown, so we apply the boundary conditions to the Hamiltonian. Note that $H(T_i) = 0$, because this is a minimum time problem with a fixed final state~\cite{FllDlvVls2012}. Additionally, $\dot{H} = 0$, because $H$ is not an explicit function of time~\cite{FllDlvVls2012}. This indicates that $H(0) = 0$, so we can solve for $\lambda(0)$. This allows us to solve for $\lambda(T_i)$ in terms of $T_i$, and we can finally solve $H(T_i)=0$ for $T_i$. This yields
\begin{align}\label{eq:proofmiet}
T_i = \frac{\sigma}{W\mathcal{W}^{-1}W^T} = \frac{\sigma_i}{d_i}.
\end{align}
Note that $T_i$ does not depend on the choice of $e_T$, so we can conclude that the minimum inter-event time is given by~\eqref{eq:proofmiet}.

\hfill $\blacksquare$
}

\emph{Proof of Theorem~\ref{th:mainResult}}

{This is a straightforward application of the general convergence result Lemma~\ref{lem:genConverge}, because Theorem~\ref{th:ldtResult} guarantees a positive minimum inter-event time so the algorithm is described by hybrid system $\mathcal{H}'$ from~\eqref{eq:modSys}.}\hfill $\blacksquare$

\begin{lemma}[Time until next event]\label{lem:timeToEvent}
For any time $t \in [t^i_\ell,t^i_{\ell+1})$, $t^i_{\ell+1}$ is next time that agent $i$ would trigger an event under~\eqref{eq:trigFunc}. That is,
\begin{align*}
t^i_{\ell+1} = \inf \{t' \geq t : \chi_i(t')=0 \text{ and } e_i(t') \neq 0\}.
\end{align*}
Under the hybrid system~\eqref{eq:system} with $\gamma_i$ defined in~\eqref{eq:clockdyn1}, the remaining time until $t^i_{\ell+1}$ is lower bounded as follows
\begin{align*}
t^i_{\ell+1} - t \geq \left\lbrace
\begin{matrix}
\frac{\sigma_i}{d_i}\left( 1 - \frac{e_i^2}{\chi_i+e_i^{2}}\right), & \text{for } (\chi_i,e_i) \neq (0,0)\\
\frac{\sigma_i}{d_i}, & \text{otherwise}
\end{matrix}
  \right. ,
\end{align*}
assuming $\chi_i \geq 0$.
\end{lemma}

\emph{Proof:}

{
This proof relies on examining each agent as a local system, as defined in~\eqref{eq:stateDyn}. We must find the minimum time to reach a point such that $\chi_i = 0$ and $e_i > 0$ from any initial point with $\chi_i(0) \geq 0$ and $e_i(0) \in \mathbb{R}$. 
We must consider three cases.
\begin{itemize}
\item \textbf{Case 1: $\chi_i(0) = 0$ and $e_i(0) = 0$}

In this case, we can use Theorem~\ref{th:ldtResult} directly, and the minimum time is $T_i$.

\item \textbf{Case 2: $\chi_i(0) > 0$ and $e_i(0) = 0$}

Note that in the proof of Theorem~\ref{th:ldtResult}, the dynamics of $\zeta$ does not depend on $\chi_i$. Therefore, with a change of coordinates $\chi_i' = \chi_i - \chi_i(0)$ we can use the same reasoning to show that the minimum time for $\chi_i$ to reach a point such that $\chi_i(t) = \chi_i(0)$ and $e_i > 0$ is $T_i$. Because $\dot{\chi}_i(0) \geq 0$ if $e_i(0) = 0$, and $\dot{\chi}_i  > 0$ if $\dot{e}_i \neq 0$, we must reach such a point before we reach a point where $\chi_i = 0$ and $e_i \neq 0$. Therefore, the minimum time is lower bounded by $T_i$.

\item \textbf{Case 3: $\chi_i > 0$ and $e_i \neq 0$}

First, we show that all points in this set are reachable using the optimal input $\mu$ from~\eqref{eq:optimalInput}. Letting $k = -\frac{e_T}{d_i T_i}$, we can write the input as
\begin{align*}
\mu = k\mathbf{1}_{|\mathcal{N}_i|},
\end{align*}
for $e_T \neq 0$. For any $k \neq 0$, then, the state will be driven along an optimal trajectory that reaches $\chi_i(T_i),e_i(T_i)= (0,e_T)$, and the value of $e_T$ will depend on $k$. Solving the differential equations with $\mu$ defined, we have
\begin{align*}
\chi_i(t) &= -(d_ik)^2 t^2 + \sigma d_ik^2 t\\
e(t) &= -d_ikt.
\end{align*}
Assuming we wish to reach the point $(\chi^*,e^*)$, we choose 
$k = \frac{e^*}{-dt}.$
This indicates that any value of $e^*\neq 0$ is reachable. Now, plugging that back into the first equation yields
\begin{align}\label{eq:chiAndTime}
\chi^* &= -e^{*2} + \frac{\sigma e^{*2}}{dt}.
\end{align}
This will indicate which values of $\chi^*$ can be reached for a given $e^*$ from the origin. Since $e^* \neq 0$, then $\lim_{t \rightarrow 0^+} \chi^* = + \infty$ and $\lim_{t \rightarrow T_i} \chi^* = e^{*2} + \frac{\sigma e^{*2}}{d\frac{\sigma}{d}}=0$. Therefore, with the proper choice of $t \in (0,T_i]$ and $k \neq 0$, our trajectory can reach any $\chi^* \geq 0$ and $e^* \neq 0$.

By the principle of optimality, the optimal input from any point in that set is still $\mu = k\mathbf{1}_{|\mathcal{N}_i|}$, and the value of $k$ will depend on the specific point.

Returning to~\eqref{eq:chiAndTime}, we solve for the time we would have reached that point on an optimal trajectory
\begin{align}\label{eq:minTime}
t &= \frac{\sigma_i e_i^2}{d_i(\chi_i+e_i^{2})}.
\end{align}
Finally, this indicates that the minimum time is
\begin{align*}
&T_i - \frac{\sigma_i e_i^2}{d_i(\chi_i+e_i^{2})}\\
= &\frac{\sigma_i}{d_i}\left( 1 - \frac{e_i^2}{\chi_i+e_i^{2}}\right).
\end{align*}

\end{itemize}

In summary, the time to reach any point such that $\chi_i = 0$ and $e_i > 0$ from any initial point with $\chi_i(0) \geq 0$ and $e_i(0) \in \mathbb{R}$ is lower bounded by
\begin{align*}
\left\lbrace
\begin{matrix}
\frac{\sigma_i}{d_i}\left( 1 - \frac{e_i^2}{\chi_i+e_i^{2}}\right), & \text{for } (\chi_i,e_i) \neq (0,0)\\
\frac{\sigma_i}{d_i}, & \text{otherwise}
\end{matrix}
  \right. .
\end{align*}

}\hfill $\blacksquare$

\emph{Proof of Theorem~\ref{th:robMIET}}

{
The lower bound on the inter-event times $\tilde{T}_i$ follows from the analysis in the proof of Lemma~\ref{lem:timeToEvent}, using the local system for agent~$i$~\eqref{eq:stateDyn}. In order for  $h_i \in [-\delta t^i, 0]$ to be satisfied  with $ \delta t^i< T_i$, we must have $e_i \neq 0$ and $\chi \geq 0$. Intuitively, the trigger condition allows events to be triggered $\delta t^i$ seconds early, so the MIET is reduced by that much. 

More formally, from the analysis in the proof of Lemma~\ref{lem:timeToEvent}, we know the time it would take to reach that point along an optimal trajectory~\eqref{eq:minTime}. Since these trajectories are optimal in a minimum time sense, there can be no faster way to reach that point from $(\chi_i,e_i)=(0,0)$, where the previous event occurred. 

For $h_i \in [-\delta t^i, 0]$, we have
\begin{align*}
\frac{\sigma_i}{d_i}\left( 1 - \frac{e_i^2}{\chi_i+e_i^{2}}\right) \in [0, \delta t^i],
\end{align*}
and so 
\begin{align*}
 \frac{\sigma_i}{d_i}\frac{e_i^2}{\chi_i+e_i^{2}} \in \left[\frac{\sigma_i}{d_i} - \delta t^i, \frac{\sigma_i}{d_i} \right].
\end{align*}
Therefore, for agent~$i$, using~\eqref{eq:minTime}, the minimum time to reach a point where the trigger condition is satisfied is
\begin{align*}
\frac{\sigma_i}{d_i} - \delta t^i.
\end{align*}

This is the MIET given in Theorem~\ref{th:robMIET}.

Now that a positive minimum inter-event time has been established, the algorithm in Theorem~\ref{th:robMIET} is known to be described by hybrid system $\mathcal{H}'$, and we can apply our general convergence result Lemma~\ref{lem:genConverge} to guarantee convergence to the set~$\{q : \hat{\phi}_i = 0\text{ } \forall\text{ } i\}$.
}\hfill $\blacksquare$


%


\begin{floatingfigure}[l]{1.25in}
\includegraphics[width=1in,height=1.25in,clip,keepaspectratio]{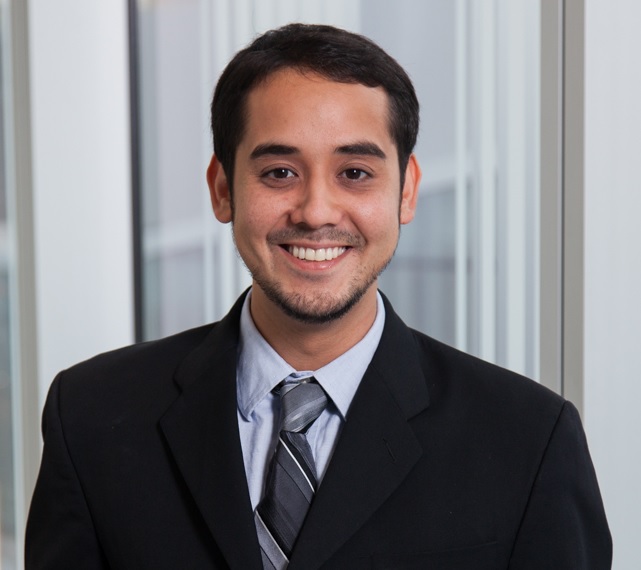}
\end{floatingfigure}
Cameron Nowzari received the Ph.D. in Mechanical Engineering from the University of California, San Diego in September 2013. He then held a postdoctoral position with the Electrical and Systems Engineering Department at the University of Pennsylvania until 2016. He is currently an Assistant Professor with the Electrical and Computer Engineering Department at George Mason University, in Fairfax, Virginia. He has received several awards including the American Automatic Control Council's O. Hugo Schuck Best Paper Award, the IEEE Control Systems Magazine Outstanding Paper Award, and the International Conference on Data Mining Best Paper Award. His current research interests include dynamical systems and control, distributed coordination algorithms, robotics, event- and self-triggered control, Markov processes, network science, spreading processes on networks, and the Internet of Things.

\hspace{-50ex}
\begin{floatingfigure}[l]{1.25in}
\includegraphics[width=1.25in,height=1in,clip,keepaspectratio]{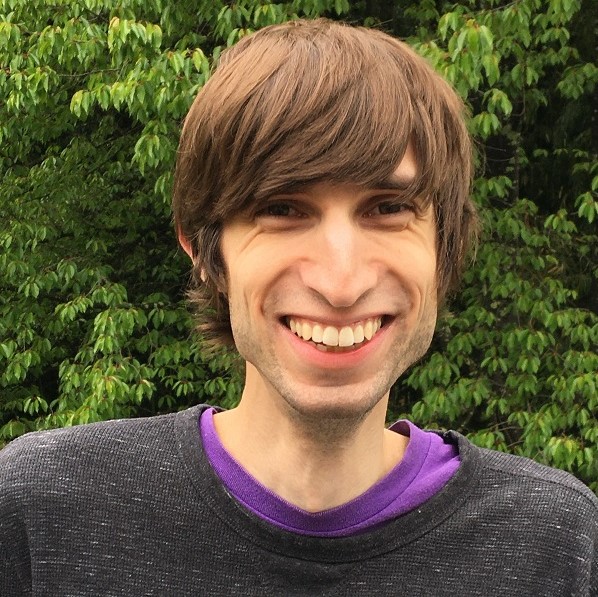}
\end{floatingfigure}
James Berneburg received the BS in Electrical Engineering from George Mason University in 2017. He is currently pursuing a Ph.D. in electrical engineering from George Mason University. He was a finalist for the Best Student Paper Award at the 2019 American Control Conference. His research interests include multi-agent systems, event-triggered control, and nonlinear control.


\end{document}